\documentclass[11pt]{article}
\usepackage{amssymb,amsfonts,amsmath,amsthm}%
\usepackage{esint}%
\usepackage{graphicx}%
\usepackage{color}%

\usepackage[font=small,format=plain,labelfont=bf,up]{caption}

\usepackage{hyperref}%
\definecolor{gray}{rgb}{0.1,0.1,.1}
\hypersetup{
    colorlinks=true,       
    linkcolor=gray,          
    citecolor=gray      
}
\newcommand{\figdraft}{false}%
\newcommand{\figfile}[1]{#1}%
\theoremstyle{plain}%
\newtheorem*{conjecture*}{Main result}%
\newcommand{\mhparagraph}[1]{\bigpar\underline{\it{#1}}.\,}

\newcommand{\Rset}{{\mathbb{R}}}

\newcommand{\ocinterval}[2]{(#1,\,#2]}%
\newcommand{\cointerval}[2]{[#1,\,#2)}%
\newcommand{\oointerval}[2]{(#1,\,#2)}%
\newcommand{\ccinterval}[2]{[#1,\,#2]}%

\newcommand{\DO}[1]{{O\at{#1}}}
\newcommand{\Do}[1]{{o\at{#1}}}

\newcommand{\evsw}{\mathrm{sw}}
\newcommand{\evsp}{\mathrm{sp}}
\newcommand{\evme}{\mathrm{me}}
\newcommand{\ev}{\mathrm{ev}}

\newcommand{\crit}{{\rm crit}}
\newcommand{\eq}{{\rm eq}}

\newcommand{\const}{{\rm const}}

\newcommand{\tdots}{{...}}%

\newlength{\mhpicDwidth}
\newlength{\mhpicDvsep}
\newlength{\mhpicDhsep}

\newlength{\mhpicPwidth}
\newlength{\mhpicPvsep}
\newlength{\mhpicPhsep}

\newlength{\mhpicWhsep}

\newcommand{\pair}[2]{{\left({#1},\,{#2}\right)}}

\newcommand{\bpair}[2]{{\big({#1},\,{#2}\big)}}
\newcommand{\Bpair}[2]{{\Big({#1},\,{#2}\Big)}}
\newcommand{\at}[1]{{\left({#1}\right)}}
\newcommand{\nat}[1]{(#1)}
\newcommand{\bat}[1]{{\big(#1\big)}}
\newcommand{\Bat}[1]{{\Big(#1\Big)}}

\newcommand{\quadruple}[4]{{\left({#1},\,{#2},\,{#3},\,{#4}\right)}}

\newcommand{\bigpar}{\par\quad\newline\noindent}

\newcommand{\jump}[1]{{|\![#1]\!|}}

\newcommand{\abs}[1]{\left|{#1}\right|}

\newcommand{\dint}[1]{\,\mathrm{d}#1}

\newcommand{\al}{{\alpha}}

\newcommand{\ga}{{\gamma}}

\newcommand{\eps}{{\varepsilon}}

\newcommand{\la}{{\lambda}}
\newcommand{\si}{{\sigma}}

\newcommand{\calD}{\mathcal{D}}
\newcommand{\calE}{\mathcal{E}}

\newcommand{\calH}{\mathcal{H}}

\newcommand{\calS}{\mathcal{S}}
\newcommand{\calT}{\mathcal{T}}

\usepackage[%
a4paper,%
nohead,%
twoside,
ignorefoot,
scale=0.8,
bindingoffset=.75cm
]{geometry}%
\usepackage{float}%
\begin{document}%
%
%
\title{ Kramers and non-Kramers Phase Transitions \\
in Many-Particle Systems with Dynamical Constraint}%
\date{\today}%
\author{%
Michael Herrmann\footnote{Universit\"at des Saarlandes, FR Mathematik,
{\tt{michael.herrmann@math.uni-sb.de}}}
\and
Barbara Niethammer\footnote{University of Oxford, Oxford Centre for Nonlinear PDE,
{\tt{niethammer@maths.ox.ac.uk}}}
\and
 Juan~J.~L.~Vel{\'{a}}zquez\footnote{Universit\"at Bonn, Institut f\"ur Angewandte Mathematik, {\tt{velazquez@iam.uni-bonn.de}}}%
}%
\maketitle
%
%
%
\begin{abstract}%
We study a Fokker-Planck equation with double-well potential that is nonlocally 
driven by a dynamical constraint and involves two small parameters.  
Relying on formal asymptotics we identify several parameter regimes
and derive reduced dynamical models for different types of phase transitions.
\end{abstract}%
%
%
\quad\newline\noindent%
\begin{minipage}[t]{0.15\textwidth}%
Keywords:%
\end{minipage}%
\begin{minipage}[t]{0.8\textwidth}%
\emph{multi-scale dynamics, gradient flows with dynamical constraint, }
\\%
\emph{phase transitions, hysteresis, Fokker-Planck equation, Kramers' formula}
\end{minipage}%
\medskip
\newline\noindent
\begin{minipage}[t]{0.15\textwidth}%
MSC (2010): %
\end{minipage}%
\begin{minipage}[t]{0.8\textwidth}%
35B40, 
35Q84, 
82C26, 
82C31 
\end{minipage}%
%
%
%
%
%
\setcounter{tocdepth}{5} %
\setcounter{secnumdepth}{3}
{\scriptsize{\tableofcontents}}%
%

\newcommand{\FPM}{FP}%
\newcommand{\TPM}{TP}%
\newcommand{\PWM}{PW}%
\newcommand{\MSM}{MS}%

%
%
\section{Introduction}\label{sec:intro}
%
%
In this paper we investigate the different dynamical regimes in a  Fokker-Planck equation with multiple scales that was introduced in \cite{DGH10} to describe the charging and discharging of lithium-ion batteries, a process that exhibits pronounced hysteretic effects \cite{DJGHMG10,DGGHJ11}. The model, to which we refer as (\FPM), governs the evolution of a statistical ensemble of identical particles and is given by the nonlocal Fokker-Plank equation
\begin{align}
\label{FPModel.PDE}\tag{\FPM$_1$}%
\tau\partial_t\varrho\pair{x}{t}&=\partial_x\Bat{\nu^2\partial_x\varrho\pair{x}{t}+\bat{H^\prime\at{x}-\si\at{t}}\varrho\pair{x}{t}}.
\end{align}
Here $H$ is the free energy of a single particle with thermodynamic state $x\in\Rset$, the probability density $\varrho\pair{\cdot}{t}$ describes the state of the whole system at time $t$, and $\si$ reflects that the system is subjected to some external forcing. Moreover, $\tau>0$ is the typical relaxation time of a single particle and  $\nu>0$ accounts for entropic effects (stochastic fluctuations).
\par
The model (\FPM) has two crucial features which cause highly nontrivial dynamics. First, the free energy $H$ is a double-well potential, hence there exist two different stable equilibria for each particle. Second, the system is not driven directly but via a time-dependent control parameter. In our case this parameter is the first moment of $\varrho$, that means we impose the dynamical constraint
\begin{align}
\int_\Rset{x}\varrho\pair{x}{t}\dint{x}&=\ell\at{t}\,,%
\label{FPModel.Constraint}\tag{\FPM$_{2}$}%
\end{align}
where $\ell$ is some given function in time, and a direct calculation shows that
\eqref{FPModel.Constraint} is equivalent to
\begin{align}
\label{FPModel.Multiplier}\tag{\FPM$_{2}^{\,\prime}$}%
\si\at{t}&=\int_{\Rset}H^\prime\at{x}\varrho\pair{x}{t}\dint{x}+\tau\dot\ell\at{t},
\end{align}
provided that the initial data satisfy $\int_\Rset{x}\varrho\pair{x}{0}\dint{x}=\ell\at{0}$.
The closure relation \eqref{FPModel.Multiplier} implies that \eqref{FPModel.PDE} is a 
nonlocal and nonlinear PDE. Well-posedness was proven in \cite{DHMRW11} on bounded domains, but we are not aware of any result about the qualitative properties of solutions. 
\par
An intriguing property of (\FPM) is that its dynamics involves
three different time scales. On the one hand, there are the relaxation time $\tau$ and the time scale of the dynamical constraint. On the other hand there is the scale of probabilistic transitions between different local minima of the effective energy $H_\si\at{x}=H\at{x}-\si{x}$, that means 
particles can move between the different wells due to stochastic fluctuations (large deviations). 
Kramers studied such transitions in the context of chemical reactions \cite{Kra40} and derived the characteristic time scale
\begin{align}
\label{Intro.KramersScale}
\tau \exp\at{\frac{\triangle{H}_\si}{\nu^2}},
\end{align}
where $\triangle H_\si$ is the minimal difference of energy between the local maximum any of the local minima of $H_\si$. In what follows we always assume that $\dot\ell$ is of order $1$, whereas both $\tau$ and $\nu$ are supposed to be small.
\par
Our goal in this paper is to identify different parameter regimes and to describe the asymptotics of (\FPM) in the limit $\nu,\tau\to0$.  To this end we focus on strictly increasing constraints and describe four different mechanisms of mass transfer between two stable regions. The corresponding
four types of phase transitions are, roughly speaking, related to two main regimes, which we refer to as
\emph{fast reaction regime} and \emph{slow reaction region}, respectively. 
The dominant effect in the fast reaction regime is mass exchange according to 
Kramers' formula. This appears for very small $\tau$ and covers, as limiting case, also
the quasi-stationary regime $\tau=0$. The slow reaction regime, however, corresponds to very small $\nu$ and Kramers' formula is not relevant anymore. Instead, phase transitions 
are dominated by transport along characteristics and this causes
rather complicated dynamics since localized peaks of mass can enter the spinodal region of $H$.
\par
In both the slow reaction and the fast reaction regimes we are able to characterize the small parameter dynamics in terms of a few averaged quantities only. Detailed descriptions of the corresponding limit models are given in the introductions to Sections \ref{sec:fast} and \ref{sec:slow}, respectively.
%
%
\subsection{Preliminaries about Fokker-Planck equations}\label{sec:prelim1}
%
Before we give a more detailed overview on the different dynamical regimes we 
specify our assumptions on $H$ and
review some basic facts about Fokker-Planck equations.
%
%
%
\subsubsection{Assumptions on the potential}
\begin{figure}[ht!]
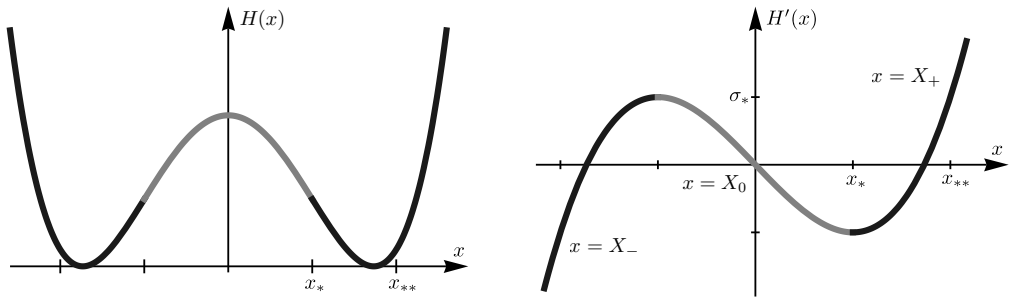
%
\centering{%
\includegraphics[width=0.4\textwidth, draft=\figdraft]%
{\figfile{potential}}%
\hspace{0.033\textwidth}%
\includegraphics[width=0.4\textwidth, draft=\figdraft]%
{\figfile{derivative}}%
}%
\caption{%
$H$ is supposed to be a 'typical' but even double well potential, as for instance \eqref{ExamplePot}.
The functions $X_-$, $X_0$, and $X_+$ denote the three branches of the inverse of $H^\prime$. 
}%
\label{Fig:Potential}%
\end{figure}

%
In this paper we assume that $H$ is an even double-well potential that
satisfies the following conditions, see Figure
\ref{Fig:Potential},
\begin{enumerate}
\item[(A1)] 
$H$ is even, sufficiently smooth (at least $C^3$), and 
$H^\prime\at{x}$ grows at least linearly as $x\to\infty$.
\item[(A2)] 
There exist constants  $0<x_*<x_{**}$ and $0<\si_*$ such that
\begin{enumerate}
\item$H^\prime\at{x_*}=-\si_*$ and $H^\prime\at{x_{**}}=\si_*$,
\item $-\si_*<H^\prime\at{x}<0$ and $H^{\prime\prime}\at{x}<0$ for $0<x<x_*$,
\item $-\si_*<H^\prime\at{x}<\infty$ and $H^{\prime\prime}\at{x}>0$ for $x>x_*$.
\end{enumerate}
In particular, the inverse of $H^\prime$ has the three strictly monotone branches 
\begin{align*}
X_-:\ocinterval{-\infty}{\si_*}\to
\ocinterval{-\infty}{-x_*}\,,\quad
X_0:\ccinterval{-\si_*}{\si_*}\to\ccinterval{-x_*}{x_*}\,,\quad
X_-:\cointerval{-\si_*}{\infty}\to
\cointerval{x_*}{\infty}\,.
\end{align*}
\item [(A3)] 
The functions $X_+\circ{H^\prime}$ and $X_-\circ H^\prime$ are concave on the spinodal interval $\oointerval{-x_*}{x_*}$.
\end{enumerate}
The assumptions in (A1) and (A2) are made for convenience and might be weakened for the price of
further technical and notational efforts; (A3) is a
 geometric condition that becomes important in the slow reaction regime that is 
discussed in Section~\ref{sec:slow}. Notice that all assumptions are in particular satisfied 
for the standard double-well potential
\begin{align}
\label{ExamplePot}
H\at{x}=\at{x^2-1}^2.
\end{align}
\par
In what follows we refer to $\oointerval{-\infty}{-x_*}$ and  
$\oointerval{x_*}{\infty}$ as the stable intervals, whereas the spinodal region
$\oointerval{-x_*}{x_*}$ is called the unstable interval. This nomenclature is motivated
by the different properties of the transport term in \eqref{FPModel.PDE}. In both stable intervals
adjacent characteristics approach each other exponentially fast, hence there is a strong tendency to 
concentrate mass into narrow peaks. In the unstable region, however, 
the separation of adjacent characteristics delocalizes any peak with positive width. 
%
%
\subsubsection{Thermodynamical aspects}
%
Fokker-Plank equations like (\FPM) are derived in \cite{DGH10} from First Principles and provide
a thermodynamically consistent model for a many-particle system with dynamical constraint. In particular, 
the second law of thermodynamics can be stated as 
\begin{align}
\notag
\dot{\calE}=-\calD+\si\dot\ell,
\end{align}
where $\calE$ is the free energy of the system and $\calD$ the dissipation. They are given
by
\begin{align*}
\calD\pair{\varrho}{\si}=\frac{1}{\tau}\int_\Rset\frac{1}{\varrho\at{x}}\Bat{\bat{H^\prime\at{x}-\si}\varrho\at{x}+\nu^2\partial_x\varrho\at{x}}^2\dint{x}\geq0
\end{align*}
and $\calE=\calH-\calS$, where
\begin{align*}
\calH\at{\varrho}=\int_{\Rset}H\at{x}\varrho\at{x}\dint{x},\qquad
\calS\at{\varrho}=-\nu^2\int_\Rset \varrho\at{x}\log\varrho\at{x}\dint{x}
\end{align*}
denote the internal energy and entropy of the many-particle system, respectively.
\par
It is well known that the Fokker-Planck equation without constraint, that is \eqref{FPModel.PDE} with $\sigma =0$, admits several interpretations as a gradient flow. There is for instance a linear structure, which
has been exploited in \cite{PSV10} in order to derive the effective dynamics in the limit 
$\nu \to 0$.  Of particular interest, however, is the nonlinear Wasserstein gradient flow structure, see \cite{JKO97,JKO98,HN11,AMPSV11}, as this structure is compatible with the constraint. 
More precisely, (\FPM) with $\dot\ell=0$ is the Wasserstein gradient flow for $\calE$ on the constraint manifold $\int_{\Rset}\varrho\at{x}\dint{x}=\ell$, and $\dot\ell\neq0$ describes a drift transversal to this manifold.
\par
The entropic term $\nu$ is often supposed to be very small but it is important that 
$\nu$ is positive. More precisely,
without the diffusive term
the qualitative properties of solutions would
strongly depend on microscopic details of the initial data, and hence it would be impossible to characterize the limit $\tau,\nu\to0$ in terms of macroscopic, i.e. averaged, quantities only, see \cite{DGH11}. A key observation is that
the singular perturbation $\nu^2\partial_x^2\varrho$  regularizes the 
macroscopic evolution, at least for some classes of initial data, in the following sense. Microscopic 
small-scale effects are still relevant on the macroscopic scale, but they are independent of the initial details and affect the system in a well-defined manner. As a consequence we now obtain a well-posed limit model
for macroscopic quantities.
Another approach to ensure well-defined macroscopic behavior is investigated in \cite{MT11}. The key idea there is to mimic entropic effects by assuming that each particle is affected by a slightly perturbed potential.
The macroscopic evolution is then completely determined by the dynamical constraint, the macroscopic initial data,
and the probability distribution of the perturbations.
%
%
\subsubsection{Dynamics in the unconstrained case}
%
We next summarize some facts about the dynamics of \eqref{FPModel.PDE} with time-independent $\si$ and small parameters $0<\nu,\tau\ll1$. For $\si\notin\ccinterval{-\si_*}{+\si_*}$, the effective potential $H_\si\at{x}=H\at{x}-\si{x}$ possesses a single critical point $\hat{x}$ that corresponds to a global minimum. The system then relaxes very fast (on the time scale $\tau$) to its unique equilibrium state 
\begin{align}
\label{Intro.Equilibrium}
\varrho_\eq\at{x}=\frac{\exp\at{-\frac{H_\si\at{x}}{\nu^2}}}{Z}\,
\end{align}
where $Z$ is a normalization constant ensuring $\int_\Rset\varrho_\eq\dint{x}=1$. This equilibrium 
density $\varrho_\eq$ has a single peak of width $\nu$ located at $\hat{x}$ and decays exponentially as $x\to\pm\infty$.
\par
The situation is different for $\si\in\oointerval{-\si_*}{+\si_*}$ since $H_\si$ now exhibits
a double well-structure with two wells (local minima) at $x_\pm=X_\pm\at{\si}$ that are separated by a barrier (local maximum) at $x_0=X_0\at\si$. Initially the system relaxes very quickly, and approaches (for smooth initial data) a state composed of two narrow peaks 
located at the wells. Both peaks have masses $m_-$ and $m_+$ with $m_-+m_+=1$, but the precise values  
of $m_\pm$ depend strongly on the initial data. This fast transition reflects that each particle in the system is strongly attracted by the nearest well due to the gradient flow structure. 
\par
The resulting state, however, is in general not an equilibrium but only a metastable state. 
The underlying  physical argument is that particles can pass the energy barrier due to stochastic fluctuations. In the generic case, in which both wells have different energies, it is of course more likely for a particle to cross the barrier coming from the well with higher energy, and thus there is a net flux of mass towards the well with lower energy.
This flux is, for small $\nu\ll1$, given by Kramers' celebrated formula and 
guarantees that the system approaches its equilibrium on the slow time scale \eqref{Intro.KramersScale}.
The corresponding equilibrium solution is again given by \eqref{Intro.Equilibrium} and has now two peaks with a definite mass distribution between the wells. Notice, however, that for small $\nu$ almost all the mass of an equilibrium solution is confined to the well with lower energy.
%
%
\subsection{Overview on different types of phase transitions}\label{sec:prelim}
%
%
Due to the different time scales, the dynamics of (\FPM) can be very complicated, and 
we are far from being able to characterize the small parameter dynamics for all types of initial data and all reasonable dynamical constraints. We thus restrict most of our considerations to strictly increasing dynamical constraints and well-prepared initial; only in Section \ref{sect:332} we allow for non-monotone constraints.
%
\subsubsection{Monotone constraints and well-prepared initial data}
%
%
%
In what follows we consider functions $\ell$ with
\begin{align}
\label{Intro.Constraint}
\ell\at{0}<-x_{**}\,,\qquad 0<c\leq\dot{\ell}\at{t}\leq{C}<\infty \quad\text{for all}\quad t\geq0,
\end{align}
where $c$, and $C$ are given constants. Since $\tau$ is small, the system then relaxes very quickly to a local equilibrium state with $\si\at{0}=H^\prime\at{\ell\at{0}}<-\si_*$. We can therefore assume that the initial mass is concentrated in a narrow peak, that means 
\begin{align}
\label{Intro.InitialData}
\varrho\pair{x}{0}\approx \delta_{\ell\at{0}}\at{x}\,,
\end{align}
where the right hand side abbreviates the Dirac distribution at $\ell\at{0}$. An even better approximation, that also accounts for the small entropic effects caused by $0<\nu\ll1$, is
\begin{align}
\label{Intro.InitialDataRefined}
\varrho\pair{x}{0}\approx\frac{1}{\nu} \sqrt{\frac{\al}{2\pi}}\exp\at{-\frac{\alpha\at{x-\ell\at{0}}^2}{2\nu^2}},\qquad
\alpha=H^{\prime\prime}\at{\ell\at0}>0.
\end{align}
The dynamical constraint implies $\dot{\si}>0$,  so that the peak starts moving to the right and the system quickly relaxes to a new local equilibrium state. For sufficiently 
small times, that means as long as $\ell\at{t}<-x_{**}$, the system can therefore be described by the \emph{single peak model} 
\begin{align}
\label{Intro.StableEvolution1}
\varrho\pair{x}{t}\approx \delta_{\ell\at{t}}\at{x},\qquad
m_-\at{t}=1,\qquad m_+\at{t}=0,\qquad\si\at{t}=H^\prime\at{\ell\at{t}}\,,
\end{align}
where we write $m_-\at{t}=1$ and $m_+\at{t}=0$ to indicate that all mass is confined in the left stable interval $\oointerval{-\infty}{-x_*}$. Moreover, at some time $t_*>0$ we have $\ell\at{t}=x_{**}$, and for $t>t_*$ the system again relaxes quickly to a local equilibrium state, that means we have
\begin{align}
\label{Intro.StableEvolution2}
\varrho\pair{x}{t}\approx \delta_{\ell\at{t}}\at{x},\qquad
m_-\at{t}=0,\qquad m_+\at{t}=1,\qquad\si\at{t}=H^\prime\at{\ell\at{t}}\,,
\end{align}
where $m_-\at{t}=0$ and $m_+\at{t}=1$ now reflect that all mass has been transferred to the second stable region $\oointerval{x_*}{\infty}$.
\bigpar
The key question is what happens between $t=0$ and $t=t_*$. 
It was already observed in \cite{DGH10} that, depending on the relation between $\tau$ and $\nu$, there are at least four types of  phase transitions
driven by rather different mechanisms of mass transfer from the left stable region into the right one. The main objective of this paper is to investigate the different regimes and to derive asymptotic formulas for the dynamics.
%
%
\subsubsection{Different regimes in numerical simulations}
%
%
\begin{figure}[ht!]
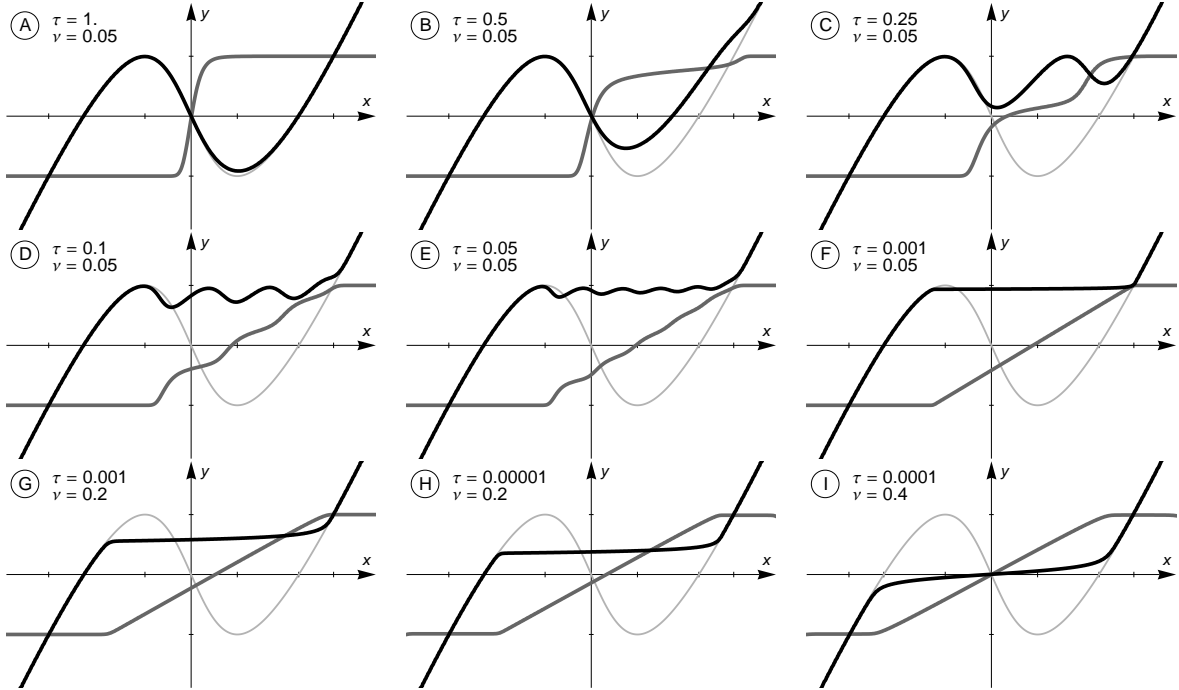
%
\centering{%
\includegraphics[width=0.3\textwidth, draft=\figdraft]%
{\figfile{ivp_sim_1}}%
\hspace{0.025\textwidth}%
\includegraphics[width=0.3\textwidth, draft=\figdraft]%
{\figfile{ivp_sim_2}}%
\hspace{0.025\textwidth}%
\includegraphics[width=0.3\textwidth, draft=\figdraft]%
{\figfile{ivp_sim_3}}%
\\%
\includegraphics[width=0.3\textwidth, draft=\figdraft]%
{\figfile{ivp_sim_4}}%
\hspace{0.025\textwidth}%
\includegraphics[width=0.3\textwidth, draft=\figdraft]%
{\figfile{ivp_sim_5}}%
\hspace{0.025\textwidth}%
\includegraphics[width=0.3\textwidth, draft=\figdraft]%
{\figfile{ivp_sim_6}}%
\\%
\includegraphics[width=0.3\textwidth, draft=\figdraft]%
{\figfile{ivp_sim_7}}%
\hspace{0.025\textwidth}%
\includegraphics[width=0.3\textwidth, draft=\figdraft]%
{\figfile{ivp_sim_8}}%
\hspace{0.025\textwidth}%
\includegraphics[width=0.3\textwidth, draft=\figdraft]%
{\figfile{ivp_sim_9}}%
}%
\caption{Numerical solutions to the initial value problem \eqref{Intro:IVP.Data} for several values of $\tau$ and $\nu$. The curves $\Gamma_{\mathrm{state}}$ and $\Gamma_{\mathrm{phase}}$ from \eqref{Intro:IVP.Results} are drawn in Black and Dark Gray, respectively; the Light Gray curve represents the graph of $H^\prime$.
}%
\label{Fig:IVP1}%
\end{figure}
\begin{figure}[ht!]
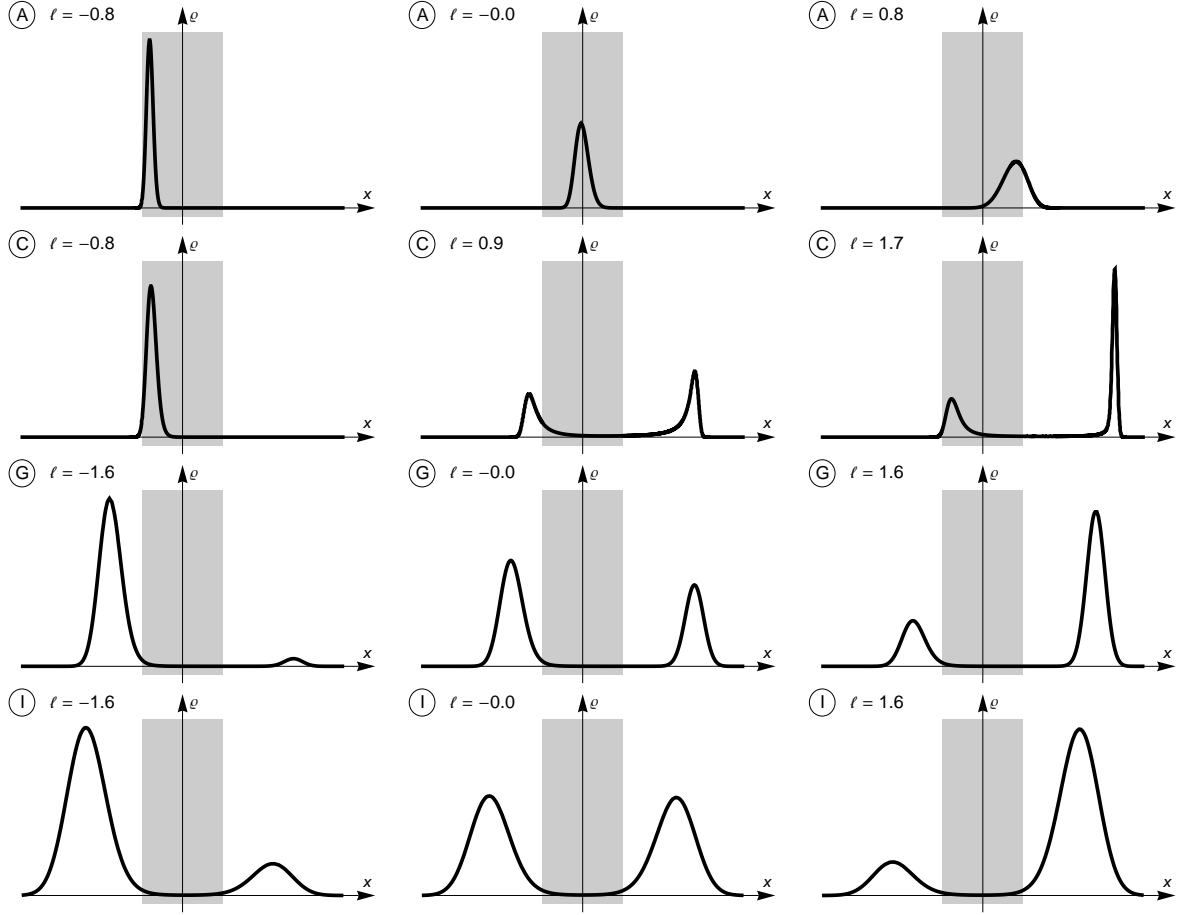
%
\centering{%
\includegraphics[width=0.3\textwidth, draft=\figdraft]%
{\figfile{ivp_dens_11}}%
\hspace{0.025\textwidth}%
\includegraphics[width=0.3\textwidth, draft=\figdraft]%
{\figfile{ivp_dens_12}}%
\hspace{0.025\textwidth}%
\includegraphics[width=0.3\textwidth, draft=\figdraft]%
{\figfile{ivp_dens_13}}%
\\%
\includegraphics[width=0.3\textwidth, draft=\figdraft]%
{\figfile{ivp_dens_21}}%
\hspace{0.025\textwidth}%
\includegraphics[width=0.3\textwidth, draft=\figdraft]%
{\figfile{ivp_dens_22}}%
\hspace{0.025\textwidth}%
\includegraphics[width=0.3\textwidth, draft=\figdraft]%
{\figfile{ivp_dens_23}}%
\\%
\includegraphics[width=0.3\textwidth, draft=\figdraft]%
{\figfile{ivp_dens_31}}%
\hspace{0.025\textwidth}%
\includegraphics[width=0.3\textwidth, draft=\figdraft]%
{\figfile{ivp_dens_32}}%
\hspace{0.025\textwidth}%
\includegraphics[width=0.3\textwidth, draft=\figdraft]%
{\figfile{ivp_dens_33}}%
\\%
\includegraphics[width=0.3\textwidth, draft=\figdraft]%
{\figfile{ivp_dens_41}}%
\hspace{0.025\textwidth}%
\includegraphics[width=0.3\textwidth, draft=\figdraft]%
{\figfile{ivp_dens_42}}%
\hspace{0.025\textwidth}%
\includegraphics[width=0.3\textwidth, draft=\figdraft]%
{\figfile{ivp_dens_43}}%
}%
\caption{Snapshots of $\varrho$ at several times corresponding to four simulations from Figure \ref{Fig:IVP1};
the shaded area represents the unstable interval.
\emph{First row}. A stable peak can pass through the unstable region.
\emph{Second row}. An unstable peak can split into two stable peaks.
\emph{Third row}. Two stable peaks exchange mass by a Kramers-type phase transition ($\si\approx\const>0$).
\emph{Fourth row}. Two stable peaks exchange mass according to the quasi-stationary dynamics ($\si\approx0$).
}%
\label{Fig:IVP2} %
\end{figure}%
We illustrate the different types of phase transitions for the driven initial value problem \eqref{Intro.Constraint} and  \eqref{Intro.InitialData} by numerical simulations with 
\begin{align}
\label{Intro:IVP.Data}
H^\prime\at{x}=x-2\arctan\at{x}\,,\qquad\ell\at{t}=t-4\,,\qquad
t\in\oointerval{0}{8}\,,
\end{align}
and initial data $\varrho\pair{\cdot}{0}$ as in \eqref{Intro.InitialDataRefined}.
Figure \ref{Fig:IVP1} visualizes the numerical solutions by means of two curves
\begin{align}
\label{Intro:IVP.Results}
\Gamma_\mathrm{state}\,:\,t\mapsto \pair{\ell\at{t}}{y\at{t}},\qquad
\Gamma_\mathrm{phase}\,:\, t\mapsto \pair{\ell\at{t}}{\si_*\bat{m_+\at{t}-m_-\at{t}}}.
\end{align}
These curves represent the \emph{macroscopic state} of the system and the (rescaled) \emph{phase field}, respectively, and are defined by
\begin{align*}
y\at{t}=\int_\Rset{H}^\prime\at{x}\varrho\pair{x}{t}\dint{x},\quad
m_+\at{t}=\int_0^{+\infty}\varrho\pair{x}{t}\dint{x},\quad
m_-\at{t}=\int_{-\infty}^0\varrho\pair{x}{t}\dint{x}.
\end{align*}
The microscopic state of the system is illustrated in Figure \ref{Fig:IVP2} by snapshots of $\varrho$ at different times.
\par
We emphasize that Figures \ref{Fig:IVP1} and \ref{Fig:IVP2} also illustrate the different dynamical regimes for fixed $0<\nu\ll1$, in the sense that the limit $\tau\to0$ can be regarded as a passage from ($A$) to ($I$). However, as our results will show, there is an exponential scale separation between the different regimes and thus it is very hard to capture all types of phase transitions in numerical simulations with the same value of $\nu$.
\par%
The numerical simulations illustrate that there exist the following types of phase transitions.
\begin{list}{$\bullet$}{%
\setlength{\leftmargin}{6\parindent}%
\setlength{\rightmargin}{\parindent}%
\setlength{\labelwidth}{3\parindent}%
\setlength{\itemindent}{0\parindent}%
\setlength{\labelsep}{1\parindent} 
}%
\item[\emph{Type $I$}, Example $(A)$:] 
At some time the narrow peak enters the unstable region due to the dynamical constraint
and starts to widen due to the separation of characteristics. However, the transport is much 
faster than the widening, so that the peak can pass trough the unstable region.
\item[\emph{Type $II$}, Examples $(B), (C), (D), (E)$:] 
The peak still enters the unstable region but now the widening is much faster than before. In particular, the unstable peak delocalizes, and the system quickly forms new peaks in each of the stable regions. At a later time the left peak enters the unstable region, and the competition between transport and widening starts once more.
\item[\emph{Type $III$}, Examples $(F), (G)$:] 
The peak does not enter the unstable region anymore. Instead, at some position in the left stable interval the peak stops moving and starts loosing mass to feed another peak in the right stable interval. This Kramers-type process happens with $\si\approx\const>0$ and
goes on until the left peak has disappeared.
\item[\emph{Type $IV$}, Example (I):] 
This is the quasi-stationary limit. The phase transition is similar to Type III but happens with $\si\approx0$. 
\end{list}
\par
Notice that the types $II$ and $III$ imply hysteresis. In fact, if we revert the situation by driving the system with $\dot\ell<0$ and $\ell\at{0}>x_{**}$, the symmetry of the problem implies that the macroscopic state and the phase field are confined to the images of $\Gamma_\mathrm{state}$ and $\Gamma_\mathrm{phase}$ under point reflections at $0$.
%
%
%
\subsubsection{Main results and organization of the paper}
%
Our main results are formulas for the macroscopic evolution
in different parameter regimes. The corresponding
scaling relations between $\tau$ and $\nu$ are
summarized in Table \ref{Tbl:ScalingRegimes} and the limit models are presented in the introductions to Sections \ref{sec:fast} and \ref{sec:slow}. 
\begin{table}[ht!]%
\centering{%
\begin{tabular}{lllcl}%
\emph{condition }&\emph{parameter}&\emph{regime} &\emph{type}&\emph{}
\\\hline\\%
$\displaystyle\tau>\frac{a_\crit}{\log1/\nu}$&&slow reactions& I &\emph{single-peak limit}
\medskip\\%
$\displaystyle\tau=\frac{a}{\log{1}/{\nu}}$&$0<a<a_\crit$&slow reactions&II&\begin{minipage}[t]{0.22\textwidth}\emph{piecewise continuous}\\\emph{two-peaks evolution}\end{minipage} 
\bigskip\\\\%
$\tau=\nu^{p}$&$0<p<\tfrac{2}{3}$&slow reactions& & OPEN PROBLEM
\bigskip\\\\%
$\tau=\nu^{p}$&$\tfrac{2}{3}<p<\infty$&fast reactions&III &\begin{minipage}[t]{0.2\textwidth}\emph{limiting case of}\\\emph{Kramers' formula}\end{minipage}
\medskip\\%
$\displaystyle\tau=\exp\at{-\frac{b}{\nu^2}}$&$0<b<h_\crit$&fast reactions&III &\emph{Kramers' formula} 
\medskip\\%
$\displaystyle\tau<\exp\at{-\frac{h_\crit}{\nu^2}}$&&fast reactions& IV&
\emph{quasi-stationary limit}
\end{tabular}%
}%
\caption{Overview on the different scaling regimes for $0<\tau,\nu\ll1$. The constant $h_\crit$ is completely determined by the properties of $H$, whereas the constant $a_\crit$ depends on both $H$ and $\ell$.}%
\label{Tbl:ScalingRegimes}%
\end{table}%
\bigpar
The rest of the paper is organized as follows. In Section~\ref{sec:fast} we derive asymptotic formulas to describe Type-$III$ transitions in the regime $\exp\at{- h_\crit/\nu^2} \ll \tau \ll \nu^{2/3}$, where $h_\crit$ denotes the energy barrier of $H$. We first  recall in Section~\ref{sec:Kramers.1} the formal derivation of Kramers' formula for the mass flux between the different wells of $H_\si$. In Section~\ref{sec:Kramers.2} we then identify the critical value for $\si=\si\pair{\tau}{\nu}$ for which such a mass flux can also take place in the constrained setting and show in Section~\ref{sec:Kramers.3} that small variations in $\si$ are sufficient to adjust the mass flux according to the dynamical constraint.  Finally, in Section~\ref{sec:Kramers.4} we discuss Type-$IV$ transitions as these can be regarded as limits of Kramers type transitions.
\par
In Section~\ref{sec:slow} we consider the scaling regime $\nu=\exp\at{-a/\tau}$ with $a>0$ and discuss Type-$II$ transitions which contain, as a special case, also Type-$I$ transitions. In Section~\ref{sect:31} we first neglect all entropic effects and introduce a simplified two-peaks model that allows to understand how two Dirac peaks
interact due to the dynamical constraint. It turns out that the dynamical constraint can stabilize a peak in the unstable region, but also that at some point a bifurcation forces both peaks to merge instantaneously. In Section~\ref{sect:32} we introduce another simplified model that accounts for the stochastic fluctuations in the unstable region and allows to understand how an unstable peak delocalizes and splits into two stable peaks. In particular, we derive an asymptotic formula for the time at which such a splitting event takes place and introduce the \emph{mass splitting problem} that determines the mass distribution between the emerging peaks. In Section~\ref{sect:33} we finally combine all result and characterize the limit dynamics as intervals of regular transport that are interrupted by several types of singular event.
%
%
%
%
\section{Fast reaction regime}\label{sec:fast}
%
In this section we show that Kramers type phase transitions are also relevant in 
presence of the dynamical constraint as long as
\begin{align*}
\exp\at{-\frac{h_\crit}{\nu^2}}\ll\tau\ll\nu^{2/3},
\end{align*}
where $h_\crit = H\at{0}-\min_{x\in\Rset}H\at{x}$ denotes the energy barrier of $H$.
The key idea is that a Kramers type phase transition occurs during a time interval $(t_1,t_2)$  in which $\si$ is positive and almost constant. During this time interval $\sigma$ only changes to order $O(\nu^2)$ but this is sufficient to accommodate the dynamical constraint. Kramers' formula therefore allows to understand phase transitions of type III,  and hence that a stable peak suddenly stops moving and starts loosing mass to feed a second stable peak. 
\par
The situation is different for $\tau\ll\exp\at{-h_\crit/\nu^2}$ since then we expect to find phase transitions of type IV, that means the mass flows towards the second well as soon as it is energetically admissible. This regime is governed by the quasi-stationary approximation but can also be regarded as a limiting case of Kramers regime.
\bigpar
Our main result concerning the fast reaction regime combines the formal asymptotics for the Kramers regime and the quasi-stationary approximation and can be stated as follows.
\begin{conjecture*}
Suppose that the dynamical constraint and the initial data satisfy \eqref{Intro.Constraint} and \eqref{Intro.InitialData}, and that $\tau$ and $\nu$ are coupled by 
\begin{align}\label{epsdef}
\tau=\exp\at{-\frac{{b}}{\nu^2}}\,
\end{align}
for some constant $b\in\oointerval{0}{h_\crit}$. Then there exists a constant $\si_b\in\oointerval{0}{\si_*}$ 
such that 
\begin{enumerate}
\item the dynamical multiplier satisfies
\begin{align*}
\si\at{t}\quad\xrightarrow{\nu\to0}\quad\left\{
\begin{array}{lclcl}
H^\prime\bat{\ell\at{t}}&&\text{for}&&t<t_1\,,\\
\si_b&&\text{for}&&t_1<t<t_2\,,\\
H^\prime\bat{\ell\at{t}}&&\text{for}&&t>t_2\,
\end{array}
\right.
\end{align*}
where $t_1$ and $t_2$ are uniquely determined by
$\ell\at{t_1}=X_-\at{\si_b}$ and $\ell\at{t_2}=X_+\at{\si_b}$,
\item the state of the system  satisfies
\begin{align*}
\varrho\pair{x}{t}\quad\xrightarrow{\nu\to0}\quad
m_-\at{t}\delta_{X_-\at{\si\at{t}}}\at{x}+m_+\at{t}\delta_{X_+\at{\si\at{t}}}\at{x}\,.
\end{align*}
where $m_+\at{t}=1-m_-\at{t}$ and
\begin{align*}
m_-\at{t}=
\left\{%
\begin{array}{lclcl}
1&&\text{for}&&t<t_1\,,\\
\displaystyle\frac{X_+\at{\si_b}-\ell\at{t}}{X_+\at{\si_b}-X_-\at{\si_b}}&&\text{for}&&t_1<t<t_2\,,\\
0&&\text{for}&&t>t_2\,.
\end{array}
\right.
\end{align*}
\end{enumerate}
Moreover, the assertions remain true
\begin{enumerate}
\item 
with $\si_b=0$ if $\tau\leq\exp\at{-\frac{h_\crit}{\nu^2}}$,
\item 
with $\si_b=\si_*$ if $\tau\ll\nu^\frac{2}{3}$ but $\tau>\exp\at{-\frac{b}{\nu^2}}$ for all $b>0$.
\end{enumerate}
\end{conjecture*}
To justify the limit dynamics we review Kramers' argument for constant $\sigma$ in Section \ref{sec:Kramers.1}.
In Section \ref{sec:Kramers.2} we then derive similar asymptotic formulas for the constrained case, which allow us to 
adjust the mass flux according to the dynamical constraints in Section \ref{sec:Kramers.3}.
Moreover, in Section \ref{sec:Kramers.4} we discuss the quasi-steady approximation, 
which governs the regime $0<\tau\ll\exp\at{-h_\crit/\nu^2}$.
\bigpar
We finally mention that the limit energy is given by
\begin{align*}
E:=m_-H\bat{X_-\at{\si}}+m_+H\bat{X_+\at{\si}},
\end{align*}
and evolves according to
\begin{align*}
\dot{E}=\si\dot\ell-D_b\chi_{\{\si=\si_b\}}\dot\ell,\qquad
D_b:=
\frac{H_{\si_b}\bat{X_-\at{\si_b}}-H_{\si_b}\bat{X_+\at{\si_b}}}{X_+\at{\si_b}-X_-\at{\si_b}}\geq0\,,
\end{align*}
where $H_\si$ is defined by $H_\si\at{x}:=H\at{x}-\si{x}$ and 
$\chi_{\{\si=\si_b\}}$ denotes the usual characteristic function.
%
%
\subsection{Kramers' formula in the unconstrained case}\label{sec:Kramers.1}
%
To derive Kramers' formula for the unconstrained case we consider the Fokker-Planck equation \eqref{TPM-1} with
fixed $\sigma \in\oointerval{-\si_*}{\si_*}$ and use the abbreviations
\begin{align*}
x_i:=X_i\at{\si},\qquad \alpha_i := \abs{H^{\prime\prime}\at{x_i}},\qquad i\in\{-,0,+\}\,.
\end{align*}
\begin{figure}[ht!]%
\centering{%
\includegraphics[width=0.4\textwidth, draft=\figdraft]%
{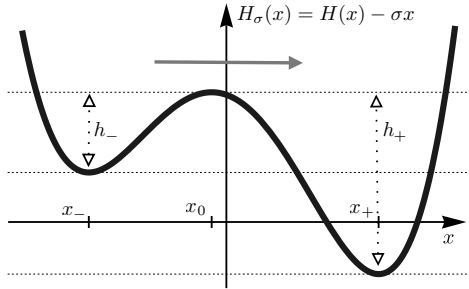}%
}%
\caption{%
Particles can cross the energy barrier between the different wells of $H_\si$ due to
random fluctuations, and this gives rise to an effective mass transfer from the well with higher energy towards 
to the well with lower energy (\emph{Kramers type phase transition}). 
}%
\label{Fig:Kramers}%
\end{figure}%
At first we approximate $\varrho$ outside the local maximum $x_0=X_0\at{\si}$ for small $\nu$ by the ansatz
\begin{align}
\label{Kramers.OuterExpansion}
\begin{split}
\varrho\pair{x}{t}\approx\left\{\begin{array}{lcl}
\mu_-\at{t}\exp\at{\displaystyle-\frac{H_\si\at{x}}{\nu^2}}&&\text{for $x<x_0$},%
\\%
\mu_+\at{t}\exp\at{\displaystyle-\frac{H_\si\at{x}}{\nu^2}}&&\text{for $x>x_0$}.%
\end{array}\right.
\end{split}
\end{align}
This is  the outer expansion and reflects the assumption that the system has a peak in either of the stable regions, where the masses are given by
\begin{align*}
m_\pm\at{t}:=\pm\mu_\pm\at{t}\int_{x_0}^{\pm\infty}
\exp\at{\displaystyle-\frac{H_\si\at{x}}{\nu^2}}\dint{x}\,.
\end{align*}
For small $\nu$ we can simplify the integrals using Laplace's method, that means we expand $H_\si$ around
$x_\pm$ to find
\begin{align*}
m_\pm\at{t}\approx
\nu\sqrt\frac{{2\pi}}{{\alpha_\pm}}\mu_\pm\at{t}\exp\at{-\frac{H_\si\at{x_\pm}}{\nu^2}}.
\end{align*}
The mass exchange between both peaks is then determined by
\begin{align}
\label{Kramers.MassTransfer}
\tau\dot{m}_\pm\at{t}=\mp R\at{t}\,,
\end{align}
where $R\at{t}$ is the mass flux at $x_0$. The key idea behind Kramers' formula is that $R\at{t}$ can be computed from the quasi-stationary approximation of $\varrho$ near $x_0$. More precisely, with the change of variables  $y=(x-x_0)/\nu$ we approximate
\begin{align*}
\nu^2\partial_x\varrho\pair{t}{x_0+\nu{y}} + H_\si^\prime\at{x_0+\nu{y}}\varrho\pair{t}{x_0+\nu{y}}\approx-R\at{t},
\end{align*}
and obtain the inner expansion
\begin{align*}
\varrho\pair{t}{x_0+\nu{y}}\approx\exp\at{-\frac{H_\si\at{x_0+\nu{y}}}{\nu^2}}\at{C\at{t}-
\frac{R\at{t}}{\nu}\int_0^y\exp\at{\frac{H_\si\at{x_0+\nu\tilde{y}}}{\nu^2}}\dint{\tilde{y}}}\,,
\end{align*}
where $C\at{t}$ is a constant of integration. For small $\nu$ and $y\gg1$ we can simplify 
the integrals by Laplace's method to obtain
\begin{align*}
\int_0^y\exp\at{\frac{H_\si\at{x_0+\nu\tilde{y}}}{\nu^2}}\dint{\tilde{y}}\approx
\int_0^\infty\exp\at{\frac{H_\si\at{x_0}-\alpha_0\tilde{y}^2}{2}}\dint{\tilde{y}}=
\sqrt{\frac{\pi}{2\alpha_0}}\exp\at{\frac{H_\si\at{x_0}}{\nu^2}},
\end{align*}
and using a similar formula for $y\ll-1$ we find
\begin{align}
\label{Kramers.InnerExpansion}\varrho\pair{t}{x_0\pm\nu\abs{y}}\approx
\exp\at{-\frac{H_\si\at{x_0\pm\nu\abs{y}}}{\nu^2}}
\at{C\at{t}\mp\frac{R(t)}{\nu}\sqrt{\frac{\pi}{2\alpha_0}}\exp\at{\frac{H_\si\at{x_0}}{\nu^2}}}\,.
\end{align}
In order to match the outer and the inner expansions, we consider $x\approx{x_0}$ and compare the asymptotic formulas \eqref{Kramers.OuterExpansion} and \eqref{Kramers.InnerExpansion}. Since both 
contain the factor $\exp\at{-H_\si\at{x}/\nu^2}$, we equate the time dependent coefficients and arrive at the
matching conditions
\begin{align*}
\mu_\pm\at{t}=C\at{t}\mp\frac{R\at{t}}{\nu}\sqrt{\frac{\pi}{2\alpha_0}}\exp\at{\frac{H_\si\at{x_0}}{\nu^2}}\,.
\end{align*}
We finally eliminate $C\at{t}$ and find the desired expression for Kramers' mass flux, namely 
\begin{align}
\label{Kramers.FluxFormula}
R\at{t}&\approx\frac{\nu\sqrt{\al_0}}{\sqrt{2\pi}}\Bat{\mu_-\at{t}-\mu_+\at{t}}\exp\at{-\frac{H_\si\at{x_0}}{\nu^2}}
\\\notag%
&\approx\frac{1}{2\pi}\exp\at{-\frac{H_\si\at{x_0}}{\nu^2}}
\at{m_-\at{t}\sqrt{\al_-\al_0}\exp\at{\frac{H_\si\at{x_-}}{\nu^2}}-m_+\at{t}\sqrt{\al_+\al_0}\exp\at{\frac{H_\si\at{x_+}}{\nu^2}}}\,.
\end{align}
Combining this with \eqref{Kramers.MassTransfer} we easily verify 
that the characteristic time for Kramers mass transfer is given by
\begin{align}
\label{Kramers.TimeScale}
\tau\exp\at{\frac{H_\si\at{x_0}-\max\at{H_\si\at{x_+},H_\si\at{x_-}}}{\nu^2}}\,.
\end{align}
We also notice that in the generic case $H_\si\at{x_+}\neq H_\si\at{x_-}$ and for small $\nu$ the mass transfer is essentially unidirectional on the time scale \eqref{Kramers.TimeScale}, that means the mass flows from the well with higher energy to the well with smaller energy, see Figure \ref{Fig:Kramers}.
%
%
\subsection{Kramers' formula in the constrained case}\label{sec:Kramers.2}
%
%
We now derive a self-consistent description for Kramers type phase transitions in the presence of the dynamical constraint. To this end it is convenient to replace $\tau$ by the parameter $b$ defined in \eqref{epsdef},
and to consider the functions
\begin{align*}
h_\pm\at\si:=H_\si\at{X_0\at\si}-H_\si\at{X_\pm\at{\si}}=
H\at{X_0\at\si}-H\at{X_\pm\at{\si}}+\si\at{X_\pm\at{\si}-X_0\at{\si}}\,.
\end{align*}
These functions are well-defined for $\abs{\si}<\si_*$ and satisfy
${\dint h_-}/{\dint\si}<0<{\dint h_+}/{\dint\si}$ with
\begin{align*}
h_-\at{-\si_*}=h_+\at{\si_*}=H\at{x_*}-H\at{x_{**}}+\si_*\at{x_*+x_{**}}=:h_*,\qquad h_+\at{-\si_*}=h_-\at{\si_*}=0\,,
\end{align*} 
and $h_-\at{0}=h_+\at{0}=H\at{0}-\min_{x\in\Rset}{H}\at{x}=:h_\crit<h_*$.
\bigpar
In terms of ${b}$ and $h_\pm$, the rescaled flux $R\at{t}/\tau$ from \eqref{Kramers.FluxFormula} can be stated as
\begin{align}
\label{Kramers.RescaledFlux}
\frac{R\at{t}}{\tau}\approx m_-\at{t} r_-\bat{\si\at{t}}-m_+\at{t} r_+\bat{\si\at{t}}\,,
\end{align}
where
\begin{align}\label{rdef}
r_\pm\at{\si}:=\frac{\sqrt{\al_\pm\at{\si}\al_0\at{\si}}}{2\pi}\exp\at{\frac{{b}-h_\pm\at{\si}}{\nu^2}}\,.
\end{align}
and $\alpha_i\at{\si}:=\abs{H^{\prime\prime}\at{X_i\at\si}}$ as above.
\bigpar%
We next present some heuristic arguments for the dynamics of the rescaled flux terms $m_-\at{t}r_-\at{t}$ and $m_+\at{t}r_+\at{t}$. To this end we assume
\begin{align*}
0<{b}<h_\crit,
\end{align*}
and recall that, due to our assumptions on the initial data, the system evolves for small times according to the single peak evolution \eqref{Intro.StableEvolution1}. In particular, the peak reaches the critical position $-x_{**}$ at time $t_0$ with $\si\at{t_0}=-\si_*$, and the dynamical constraint implies $\dot{\si}\at{t_0}>0$. Assuming that $\si$ changes regularly at $t_0$, we then conclude that
\begin{align*}
h_-\at{\si\at{t_0}}=h_*>0,\quad\frac{\dint}{\dint{t}}{h_-\at{\si\at{t}}}<0,\quad
h_+\at{\si\at{t_0}}=0,\quad\frac{\dint}{\dint{t}}{h_+\at{\si\at{t}}}>0.
\end{align*}
Consequently, for small $\nu$ and sufficiently small times $t>t_0$ we expect to find $r_-\at{t}\ll1\ll r_+\at{t}$, so crossing the energy barrier is very likely for a particle in the right well but very unlikely for a particle in the left well. However, the net transfer across the energy barrier is very small since there are essentially no particles in the right well. We thus expect that the partial masses $m_\pm$ stay constant in the limit $\nu\to0$, so the system can still be described by the single peak approximation \eqref{Intro.StableEvolution1}. At some later time $\tilde{t}_0>t_0$ 
we have $\si\at{\tilde{t}_0}=0$ and the fluxes $r_-\nat{\tilde{t}_0}$ and $r_+\nat{\tilde{t}_0}$ have the same order of magnitude. However, both are very small due to ${b}<h_\crit=h_\pm\at{\tilde{t}_0}$,  and so there is, for small $\nu$, still no effective mass transfer between the two wells of $H_\si$. 
\par
The situation changes completely at time $t_1$ defined by $\si\at{h_-\at{t_1}}={b}$. At this time, $r_-\at{t_1}$ becomes suddenly of order one and we can no longer neglect particles that move from the left well to the right one. The other flux $r_+$, however, is now very small as it is very unlikely that a particle moves the other way around.
\bigpar
As explained above, the main idea in the dynamical case is that there exist a constant $\si_b$ and a time $t_2$ such that $\si\at{t}\approx\si_{b}$ for all $t_1<t<t_2$. Kramers' mass flux can hence stabilize to continuously transfer mass from the left well to the right one. At time $t_2$, all mass has been transferred to the right well and the system again evolves according to the single peak evolution, now given by \eqref{Intro.StableEvolution2}.
\par
Before we describe the details of the mass transfer we proceed with two remarks. First, the above assumption ${b}<h_\crit$ is truly necessary: For $h_\crit<{b}<h_*$ there is still a time $t_1$ with $h_-\at{\si\at{t_1}}={b}$, but then we have $h_+\at{\si\at{t_1}}<h_-\at{\si\at{t_1}}$ and hence $r_+\at{t_1}\gg r_-\at{t_1}$, which shows that a net transfer from the left well to the right one is impossible. Moreover, for ${b}>h_*$ both 
$r_-$ and $r_+$ are always very large. In both cases we expect that the phase transition
occurs when $\si\approx0$ and is not governed by Kramers' formula anymore but by the quasi-stationary approximation. 
\par
Second, the mass flux is already determined by the dynamical constraint and the assumption $\si\at{t}\approx\si_{b}$. In fact, the constraint \eqref{FPModel.Constraint} implies that 
\begin{align*}
\ell\at{t}\approx X_-\at{\si_{b}}{m}_-\at{t} +  X_+\at{\si_{b}}{m}_+\at{t}
\end{align*}
and thus, since $\si$ does not change much, 
\begin{align*}
\dot\ell\at{t}\approx X_-\at{\si_{b}}\dot{m}_-\at{t} +  X_+\at{\si_{b}}\dot{m}_+\at{t}=
\Bat{X_-\at{\si_{b}}-X_+\at{\si_{b}}}\dot{m}_-\at{t}\,.
\end{align*}
As a consequence we obtain 
\begin{align}\label{mdot}
\dot{m}_-\at{t}=-\frac{\dot\ell\at{t}}{X_+\at{\si_{b}}-X_-\at{\si_{b}}}
\end{align}
for $t_2<t<t_3$, where $t_3>t_2$ is defined by $m_-\at{t_3}=0$, and using $\ell\at{t_2}\approx X_-\at{\si_{b}}$ it is easy to check that $\ell\at{t_3}\approx X_+\at{\si_{b}}$. 
%
%
\subsection{Adjusting the mass flux by small variations of the multiplier}\label{sec:Kramers.3}
%
%
To identify the formulas that relate the mass flux self-consistently to small temporal changes in $\si$ we consider only times $t$ with $t_1<t<t_2$ and assume that both
$m_-$ and $m_+$ are strictly positive. Of course, in order to match the resulting approximations for $\varrho$ to the single-peak evolution for $t<t_1$ and $t>t_2$ we must introduce transition layers at $t\approx{t_1}$ and $t\approx{t_2}$ corresponding to $m_-\approx1$ and $m_-\approx0$, respectively, but since these transition layers do not contribute to the limit model, we do not investigate them in detail.
\mhparagraph{Case 1 : $ \lim_{\nu \to 0} {b} >0$} %
The critical value $\si_b$ is defined by ${b}=h_-\at{\si_{{b}}}$. Thanks to ${b}<h_\crit$ we find
$0<\si_{b}<{\si_*}$ and $h_-\at{\si_{{b}}}>h_\crit>h_+\at{\si_{{b}}}$, 
so \eqref{rdef} yields
\begin{align}
\label{KramersEffectiveFluxes}
r_-\at{t}\approx\frac{\sqrt{\al_-\at{\si_{b}}\al_0\at{\si_{b}}}}{2\pi}\exp\at{\frac{h_-\at{\si_{b}}-h_-\at{\si\at{t}}}{\nu^2}},\qquad
r_+\at{t}\approx0.
\end{align}
Since $r_-\at{t}$ must be of order one, we introduce the rescaled multiplier 
\begin{align*}
\psi\at{t}:=\frac{\si\at{t}-\si_{b}}{\nu^2},
\end{align*}
and simplify \eqref{KramersEffectiveFluxes} by expanding $h_-\at{\si}$ around $\si_{b}$.
Using \eqref{Kramers.MassTransfer} and \eqref{Kramers.RescaledFlux} we then conclude that
Kramers' formula implies the mass transfer law 
\begin{align*}
-\dot{m}_-\at{t}\approx m_-\at{t}\frac{\sqrt{\al_-\at{\si_{b}}\al_0\at{\si_{b}}}}{2\pi}\exp\at{\abs{h_-^\prime\at{\si_{b}}}\psi\at{t}}\,.
\end{align*}
Comparing this with \eqref{mdot} we finally conclude that Kramers type phase transitions comply with the dynamical constraint if and only if 
\begin{align}\label{psiformula}
 \psi\at{t}\approx
\frac{\ln\at{\displaystyle\frac{2\pi}{\sqrt{\al_-\at{\si_{b}}\al_0\at{\si_{b}}}
\bat{X_+\at{\si_{b}}-X_-\at{\si_{b}}}}\frac{\dot{\ell}\at{t}}{m_-\at{t}}}}{\abs{h_-^\prime\at{\si_{b}}}}\,.
\end{align}
This is the heart of our argument. If $\psi$ evolves according to \eqref{psiformula}, then $\si-\si_{b}$ is of order $\nu^2$, and Kramers' formula 
provides a mass flux that satisfies the dynamical constraint.
\mhparagraph{Case 2 : $ \lim_{\nu \to 0} {b} =0$} %
In this limiting case, the phase transition happens when $\sigma$ is close to $\sigma_*$,  so
both $x_-$ and $x_0$ are close to $-x_*$. Moreover, the constants $\alpha_-$ and $\alpha_0$ approach zero, and hence we can no longer use the asymptotic expressions from the first case. However,
if we expand all relevant quantities around $\sigma_*$ and $-x_*$, it is still possible to derive an asymptotic formula for Kramers' flux that is consistent with the dynamical constraint.
\par
Thanks to the identities $X_-\at{\si_*}=X_0\at{\si_*}=-x_*$, $H^{\prime\prime}\at{-x_*}=0$ 
and $\ga:=-H^{\prime\prime\prime}\at{-x_*}>0$, we deduce from the definition of $X_-$ and $X_0$ that 
\begin{align*}
X_-\at{\si}=-x_*-\sqrt{\frac{2\at{\si_*-\si}}{\ga}}+\DO{\si_*-\si},\qquad
X_0\at{\si}=-x_*+\sqrt{\frac{2\at{\si_*-\si}}{\ga}}+\DO{\si_*-\si}\,.
\end{align*}
To leading order in $\si_*-\si$, we therefore find
\begin{align*}
\sqrt{\alpha_-\at{\si}\alpha_+\at\si}\approx\sqrt{2\ga\at{\si_*-\si}}\,,
\end{align*}
as well as
\begin{align*}
H_\si\bat{X_0\at\si}-H_\si\bat{X_-\at{\si}}&\approx-\frac{\ga}{6}\Bat{\bat{X_0\at{\si}+x_*}^3-
\bat{X_-\at{\si}+x_*}^3}+\at{\si_*-\si}\bat{X_0\at{\si}-X_+\at{\si}}
\\%
&=\frac{1}{3}\at{\si_*-\si}\bat{X_0\at{\si}-X_+\at{\si}}\approx\frac{4\sqrt{2}}{3\sqrt{\ga}}\at{\si_*-\si}^{3/2},
\end{align*}
so \eqref{Kramers.FluxFormula} can be simplified to
\begin{align}
\label{Kramers.FluxAsymptotics}
R\at{t}\approx\frac{m_-\at{t}\sqrt{\ga}}{\sqrt\pi}\bat{\si_*-\si\at{t}}^{1/2}\exp\at{-\frac{4\sqrt{2}}{3\sqrt{\ga}}\frac{\bat{\si_*-\si\at{t}}^{3/2}}{\nu^2}}\,.
\end{align}
In Kramers' regime this flux $R\at{t}$ should be of order $\tau$. On the other hand,
the asymptotic formula \eqref{Kramers.FluxFormula} holds only if $H_\si\at{x_0}-H_\si\at{x_-}\gg\nu^2$, and thus we shall guarantee that $\nu^{-2}\at{\si_*-\si}^{3/2}\gg1$. Both conditions can be satisfied if 
\begin{align*}
\tau\ll\nu^{2/3}\,.
\end{align*}
In fact, if we define for given $\nu\ll1$ and $\tau\ll\nu^{2/3}$ the large parameter $K$ 
by
\begin{align*}
\frac{\sqrt{\ga}}{\sqrt\pi}{K}\exp\at{-\frac{4\sqrt{2}}{3\sqrt{\ga}}K^3}=\tau\nu^{-2/3}
\end{align*}
then $R\at{t}$ becomes of order $\tau$ if $\si_*-\si$ is of order $K^2{\nu^{4/3}}\sim \at{\nu^2\ln\at{1/\nu}}^{2/3}$.
\par
Finally, we proceed as in the first case. Inserting \eqref{Kramers.FluxAsymptotics} into \eqref{Kramers.MassTransfer} gives a formula for $\dot{m}_-\at{t}/m_-\at{t}$ in terms of $\si_*-\si\at{t}$. On the other hand, the 
dynamical constraint implies 
\begin{align*}
\dot{m}_-\at{t}/m_-\at{t}\approx\at{ -\frac{\dot{\ell}\at{t}}{2x_*}}\,/\,
\at{1-\frac{\ell\at{t}+x_*}{2x_*}}
\end{align*}
via \eqref{mdot}, and eliminating $\dot{m}_-\at{t}/m_-\at{t}$ we obtain
$\si_*-\si\at{t}=\DO{K^2\nu^4/3}$ in terms of $\ell\at{t}$.
%
%
\subsection{Phase transitions in the quasi-stationary limit  }\label{sec:Kramers.4}
%
%
To conclude this section we show that the quasi-stationary approximation of (\FPM) describes 
phase transitions with $\si\at{t}\approx0$. Notice that such Type-IV transitions can be regarded
as limits of Type-III transitions in the sense that  $\si_b\searrow0$ as $b\nearrow h_\crit$.
\par
In the quasi-stationary limit we approximate $\varrho$ by the equilibrium solution \eqref{Intro.Equilibrium}  that corresponds to the current value of $\si$. In other words, we set 
\begin{align*}
\varrho\pair{x}{t}\approx\varrho_{\si\at{t}}\at{x},\qquad
\varrho_\si\at{x}:=\frac{\exp\at{\displaystyle-\frac{H\at{x}+\si{x}}{\nu^2}}}{Z\at{\si}}\,,
\end{align*}
where the normalization factor $Z\at{\si}:=Z_-\at{\si}+Z_+\at{\si}$ is given by
\begin{align}
\label{Kramers.QSA.interals}
Z_\pm\at{\si}:=\mp\int_{\pm\infty}^{X_0\at\si}\exp\at{\displaystyle-\frac{H\at{x}+\si{x}}{\nu^2}}\dint{x}\,.
\end{align}
The dynamical multiplier $\si\at{t}$ is then determined by the dynamical constraint via
\begin{align*}
\ell\at{t}=\int_\Rset\varrho_{\si\at{t}}\at{x}\dint{x}\,.
\end{align*}
We now derive asymptotic formulas that characterize the quasi-stationary dynamics for small $\nu$. At first we notice
that Laplace's method applied to \eqref{Kramers.QSA.interals} with $\si\at{t}\neq0$  yields
\begin{align*}
Z_\pm\at{\si}\approx\frac{\nu\sqrt{2\pi}}{\sqrt{\alpha_\pm\at\si}}\exp\at{\frac{-H\bat{X_\pm\at{\si}}+\si{X}_\pm\at{\si}}{\nu^2}}\,.
\end{align*}
and hence
\begin{align*}
\varrho_\si\at{x}\quad\xrightarrow{\nu\to0}\quad\left\{\begin{array}{lcl}%
\delta_{X_-\at{\si}}\at{x}&&\text{for $\si<0$},\\
\delta_{X_+\at{\si}}\at{x}&&\text{for $\si<0$}.
\end{array}\right.
\end{align*}
In particular, the system evolves according to the single-peak approximation as long as $\si$ has a sign, and 
a phase transition can occur only for $\si=0$. To describe the details of such a transition it is convenient to
rescale $\si$ by  $\si=\nu^2\psi$. 
\par
Using the expansion
\begin{align*}
H\at{X_\pm\at{\si}+\nu y}-\si\at{X_\pm\at{\si}+\nu{y}}=H\bat{X_\pm\at\si}
+\tfrac{1}{2}H^{\prime\prime}\bat{X_\pm\at{0}}\nu^2y^2-\si{X}_\pm\at{0}+
\DO{\si^2,\si\nu,\nu^3}
\end{align*}
and employing Laplace's method once more, we find
\begin{align*}
Z_\pm\at{\si\at{t}}\approx\frac{\nu\sqrt{2\pi}}{\sqrt{\alpha_\pm\at0}}\exp\at{\frac{-H\bat{X_{\pm}\at{0}}}{\nu^2}+{X}_\pm\at{0}\psi\at{t}}\,,
\end{align*}
and hence
\begin{align*}
m_-\at{t}/m_+\at{t}\approx
Z_-\at{\nu^2\psi\at{t}}/Z_+\at{\nu^2\psi\at{t}}\approx\exp\bat{-2X_+\at{0}\psi\at{t}}\,.
\end{align*} 
On the other hand, the dynamical constraint, see \eqref{mdot}, provides
\begin{align*}
m_-\at{t}\approx \frac{X_+\at{0}-\ell\at{t}}{2X_+\at{0}},\qquad
m_+\at{t}\approx \frac{X_+\at{0}+\ell\at{t}}{2X_+\at{0}},
\end{align*}
and we conclude that the rescaled multiplier evolves according to
\begin{align*}
\psi\at{t}=-\frac{\ln\at{\displaystyle\frac{X_+\at{0}-\ell\at{t}}{X_+\at{0}+\ell\at{t}}}}{2X_+\at{0}}.
\end{align*}
Notice that this formula is well defined for all times $t$ with $t_1<t<t_2$, where
$t_1$ and $t_2$ are defined by $\ell\at{t_1}=X_-\at{0}=-X_+\at{0}$ and 
$\ell\at{t_2}=+X_+\at{0}$, and satisfy $\lim_{t\to t_1}\psi\at{t}=-\infty$ and 
$\lim_{t\to t_2}\psi\at{t}=+\infty$.
%
%
%
\section{Slow reaction regime}\label{sec:slow}
%
%
This section concerns the effective dynamics of (\FPM) in the slow reaction regime: Both $\tau$ and $\nu$ are still supposed to be small, but $\nu$ is so small that `reactions', that means continuous mass transfer between the stable regions as described by Kramers' formula, are not relevant anymore. Instead, the dominant effect in (\FPM) is now transport along characteristics and this gives rise to new phenomena. In particular, localized peaks can enter the unstable region and peaks can split or merge rapidly. Notice, however, that the small entropic effects caused by $\nu>0$ are still relevant and cannot be neglected. They guarantee that each peak entering the unstable region is basically a rescaled Gaussian, and hence that such a peak behaves in a well-defined manner.
\par
We now introduce an informal concepts that is motivated by numerical simulations and turns out
to be useful for describing the asymptotic dynamics in the slow reaction regime.
We say the system is in a \emph{two-peaks configuration} if there exist two positions
$x_1,x_2$ and two masses $m_1,m_2$ with $0\leq{m_2}=1-m_1\leq1$ such that the state 
can be approximated by
\begin{align}
\notag
\varrho\pair{x}{t}\approx m_1\at{t}\delta_{x_1\at{t}}\at{x}+ m_2\at{t}\delta_{x_1\at{t}}\at{x}\,.
\end{align}
A two-peaks configuration is called \emph{stable-stable} if $x_1<-x_*$ and $x_2>x_*$, but \emph{unstable-stable} if $-x_*<x_1<x_*$ and $x_2>x_*$. Moreover, in case that one of the masses vanishes, we refer to a two-peaks configuration as a (\emph{stable} or \emph{unstable}) \emph{single-peak configuration}.
\begin{figure}[ht!]%
\centering{%
\includegraphics[width=0.975\textwidth, draft=\figdraft]%
{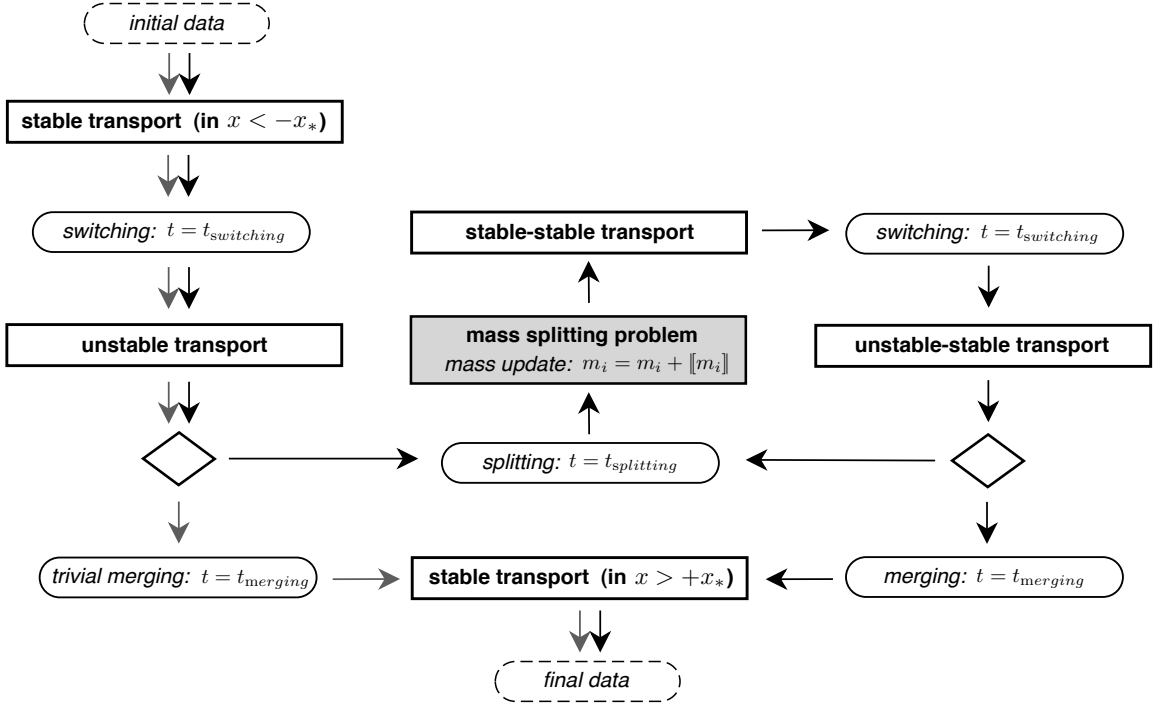}%
}%
\caption{%
Flowchart for the limit dynamics in the slow-reaction regime with strictly increasing dynamical constraint. Intervals of quasi-stationary transport (of either a single-peak or a two-peaks configuration) are interrupted by several types of singular times (corresponding to switching, merging and splitting of peaks). Gray and Black arrows indicate Type-I and Type-II phase transitions, respectively.
}%
\label{Fig:flow_chart}%
\end{figure}%
\bigpar
Our main result in this section is an asymptotic description of the
dynamics in a certain parameter regime for $\tau$ and $\nu$. The corresponding limit model is illustrated in Figure \ref{Fig:flow_chart} and can be summarized as follows.
\begin{conjecture*}
Suppose that the dynamical constraint and the initial data satisfy \eqref{Intro.Constraint} and \eqref{Intro.InitialData}, and that $\tau$ and $\nu$ are coupled by 
\begin{align}
\label{Slow:Scaling}
\nu=\exp\at{-\frac{a}{\tau}}
\end{align}
for some constant $a>0$. Then, in the limit $\tau\to0$ the dynamics can be described in terms of single-peak and two-peaks configurations according to one of the following scenarios.
\begin{enumerate}
\item 
\emph{Transport with splitting (Type-II phase transition)}:
There are intervals of quasi-stationary transport which are interrupted by singular times.
More precisely, depending on $\ell$  and the value of the parameter $a$ there exists an integer $K\geq1$ such that 
there are $K+1$ \emph{switching times}, $K$ \emph{splitting times}, and a final \emph{merging time}.
\begin{enumerate}
\item During the intervals of quasi-stationary transport, the peaks do not exchange mass and are move according to the dynamical constraint.
\item 
At each switching time a stable peak reaches the position $-x_*$ to become unstable; each switching time is followed by a splitting or the merging time.
\item 
At  each splitting time an unstable peaks splits and its mass is instantaneously transferred to the stable regions. In particular, the system jumps from an unstable-stable two-peaks configuration (or, initially, from an unstable single-peak configuration) to an emerging stable-stable two-peaks configuration. The precise values for the 
masses and the positions of the emerging peaks are determined by a \emph{mass splitting problem}.
\item 
At the final merging time the two peaks of an unstable-stable configuration
merge -- either continuously or discontinuously -- to form a single stable peak located in the region $x\geq{x_*}$. 
\end{enumerate}
\item 
\emph{Pure transport (Type-I phase transition)}:
The system is always in a single peak configuration and there exist only two singular times corresponding to switching (the peaks enters the unstable interval) and trivial continuous merging (the peak leaves the unstable interval).
\end{enumerate}
\end{conjecture*}
Notice that  it is practically impossible to perform numerical simulations with $\nu$ exponentially small in $\tau$ and 
$\tau$ much smaller than $1$. The phase transitions of type I and type II 
presented in Figures \ref{Fig:IVP1} and \ref{Fig:IVP2} are therefore not directly covered by our limit model, but a close look to the numerical data reveals that they are likewise dominated by  the interplay between
transport and widening of unstable peaks. It remains a challenging task to derive next order 
corrections to replace the singular times by transition layers with width depending on $\tau$. Another interesting but 
open question is whether there exist scaling laws different from \eqref{Slow:Scaling} that give rise to other reasonable slow reaction limits.
\bigpar
To justify the limit dynamics we introduce several reduced models that allow us to study each of the different phenomena in a simplified setting. In Section \ref{sect:31} we derive a \emph{two-peaks model} that describes the transport of two peaks in terms of a simple ODE system. Moreover, this model also reveals that the separated peaks in an unstable-stable configuration can merge instantaneously due to the dynamical constraint. In Section \ref{sect:32} we investigate the splitting of unstable peaks. To this end we propose  a \emph{peak-widening model} that accounts for the small entropic effects and allows us to derive a deterministic equation for the width of an unstable peak. In particular, we show that for small $\tau$ it can happen that this width blows up almost instantaneously, which in turn gives rise to rapid mass transfer from the unstable towards the stable regions. We then discuss the \emph{mass splitting problem}, which consists of solving a nonlocal transport equation in order to determine how much mass is transferred to each of the stable regions. Finally, in Section \ref{sect:33} we combine all partial results and characterize the limit dynamics of the original model (\FPM). In particular, we derive explicit formulas for the iterative computation of all switching, splitting, and merging times. Moreover, at the end we sketch the slow-reaction limit for non-monotone dynamical constraints.
%
%
\subsection{Two-peaks approximation: Transport and merging of peaks}\label{sect:31}
%
%
A major tool in our analysis is a simple \emph{two-peaks model} (\TPM) which describes the essential dynamics of (\FPM) as long as the system is a two-peaks configuration. The model governs the evolution of the peak positions $x_1$ and $x_2$, and reads
\begin{align}
\label{TPM-1}\tag{\TPM$_1$}
\tau\dot{x}_{1/2}&=\si-H^\prime\at{x_{1/2}},\\
\label{TPM-2a}\tag{\TPM$_2$}
\ell&=m_1x_1+m_2x_2\\
\label{TPM-2b}\tag{\TPM$_2^{\,\prime}$}
\si&=m_1H^\prime\at{x_1}+m_2H^\prime\at{x_2}+\tau\dot\ell\,.
\end{align}
Here $m_1$ and $m_2=1-m_1$ denote the constant masses of the peaks. Notice that \eqref{TPM-1} is just the characteristic ODE for \eqref{FPModel.PDE} with $\nu=0$, and this implies that (\TPM) is also a constrained gradient flow corresponding to the energy
\begin{align*}
E:=m_1H\at{x_1}+m_2H\at{x_2}\,.
\end{align*} 
In particular, using the dissipation 
\begin{align*}
D:=\tau{m_1}\dot{x}_1^2+\tau{m_2}\dot{x}_2^2=\frac{m_1}{\tau}\bat{H^\prime\at{x_1}-\si}^2+\frac{m_2}{\tau}\bat{H^\prime\at{x_2}-\si}^2,
\end{align*} 
the energy balance is given by $\dot{E}=-D+\si\dot\ell$.
\bigpar%
Since $\tau$ is small it seems natural to neglect the time derivatives in \eqref{TPM-1} and \eqref{TPM-2b}. This gives rise to the \emph{quasi-stationary approximation} to (\TPM), which consists of the algebraic equations
\begin{align}
\label{TPM.QS}
H^\prime\at{x_1}=H^\prime\at{x_2}=\si,
\qquad m_1x_1+m_2x_2=\ell.
\end{align}
For our analysis it is important to understand in which sense \eqref{TPM.QS} approximates (\TPM).
This problem is not trivial because the non-invertibility of $H^\prime$ implies that
\eqref{TPM.QS} has multiple solutions for $\abs{\ell}<x_{**}$, and thus we have to understand 
which solution branches are dynamically selected by solutions to (\TPM). 
\bigpar
In the limit $\tau\to0$, solutions to (\TPM) exhibit two important dynamical phenomena which 
correspond to changing the solution branch of \eqref{TPM.QS}. Both phenomena seem  
to be counter-intuitive at a first glance but are a consequence of the dynamical constraint. 
They can be described as follows.
\par
First, the dynamical constraint can drive the system from a stable-stable configuration to an unstable-stable configuration, that means the stable peak at $x_2$ can stabilize an unstable one at $x_1$ due to the dynamical constraint. The second effect is that this stabilization can break down eventually. When this happens, the separated peaks merge almost instantaneously to form a single stable peak.
\par%
To describe both phenomena we consider times $t\geq{t_0}$ and suppose that the dynamical constraint $\ell\at{t}$ is smooth and strictly increasing. We also suppose that the initial data for (\TPM) are well prepared via
\begin{align}
\label{TPM.IV}
x_1\at{t_0}=X_-\at{\si\at{t_0}},\quad 
x_2\at{t_0}=X_+\at{\si\at{t_0}},\quad
\ell\at{t_0}=m_1x_1\at{t_0}+m_2x_2\at{t_0}\,.
\end{align}
To elucidate the key ideas we proceed with discussing some numerical results and present some semi-rigorous analytical considerations afterwards.
%
%
%
\subsubsection{Dynamics of the two-peaks model}
%
%
%
\begin{figure}[t!]
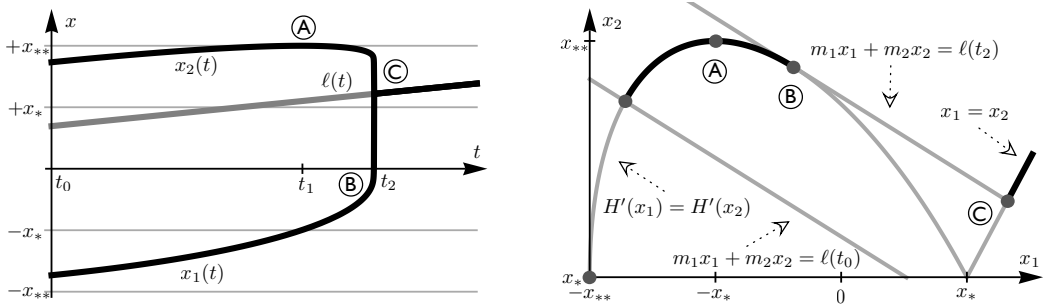
%
\centering{%
\includegraphics[width=0.4\textwidth, draft=\figdraft]%
{\figfile{two_particles}}%
\hspace{0.05\textwidth}%
\includegraphics[width=0.4\textwidth, draft=\figdraft]%
{\figfile{steady_states}}%
}%
\caption{%
\emph{Left panel:} Numerical solution to the two-peaks model (\TPM) with $0<\tau\ll1$. As long as $t-t_0$ is sufficiently small, each peak is located in one of the stable regions, but at time $t_1$ the peak at $x_1$ enters the unstable region (switching configuration $A$). At a later time $t_2$ the quasi-stationary two-peaks approximation ceases to exist, and the system jumps almost instantaneously from configuration $B$ to the single-peak configuration $C$ (discontinuous merging of peaks). \emph{Right panel:} Cartoon of all steady states in the $\pair{x_1}{x_2}$-plane. The solution from the left panel is initially confined to the curve $H^\prime\at{x_1}=H^\prime\at{x_2}$ with $x_2>x_1$ but jumps to the diagonal 
$x_1=x_2$ at time $t_2$. 
}%
\label{Fig:HyperbolicInstability}%
\end{figure}%
The left panel of Figure \ref{Fig:HyperbolicInstability} depicts a typical numerical solution to (\TPM) with $0<\tau\ll1$, potential \eqref{ExamplePot}, and initial data as in \eqref{TPM.IV}. The simulation reveals the existence of two critical times $t_1$ and $t_2$ that separate the three different regimes
\begin{align*}
\begin{array}{ccrcccrcrcccr}
R_1&\quad:\quad&t_0&<&t&<&t_1,&\;&-x_{**}&<&x_1&<&-x_{*},\\
R_2&\quad:\quad&t_1&<&t&<&t_2,&\;&-x_*&<&x_1&<&+x_{*},\\
R_3&\quad:\quad&t_2&<&t&<&t_3,&\;&+x_*&<&x_1&<&+x_{**},
\end{array}
\end{align*}
where $t_3$ denotes the final simulation time. For all times $0<t<t_3$ we have
\begin{align*}
H^\prime\at{x_1}\approx\si,\qquad H^\prime\at{x_2}\approx\si,\qquad
x_*<x_2<x_{**},
\end{align*}
and thus we expect that in the limit $\tau\to0$ the system can in fact be described by quasi-stationary peaks. The details, however, 
are different for $R_1$, $R_2$, and $R_3$. More precisely, initially we have
\begin{align*}
R_1\quad:\quad x_1\at{t}\approx X_-\at{\si\at{t}},\qquad
x_2\at{t}\approx X_+\at{\si\at{t}},
\end{align*}
that means the solution resembles a stable-stable configuration, and 
$\dot\ell>0$ implies $\dot{x}_1>0$, $\dot{x}_2>0$, and $\dot{\si}>0$.
At time $t_1$, which is defined by $x_1\at{t_1}=-x_*$, the peak at $x_1$ enters the unstable region and the configuration switches to unstable-stable. This means 
\begin{align*}
R_2\quad:\quad x_1\at{t}\approx X_0\at{\si\at{t}},\qquad x_2\approx X_+\at{\si\at{t}},
\end{align*} 
and hence
$\dot{x}_1>0$, $\dot{x}_2<0$, and $\dot{\si}<0$. 
\par%
At the second critical time $t_2$, the two-peaks approximation breaks down, that 
means the system can no longer be described by two separated peaks.
Instead, both peaks merge discontinuously, in the sense that the system jumps
from an unstable-stable two-peaks configuration to a stable single-peak configuration. The evolution for $t>t_2$ is still quasi-stationary, but involves only a single peak that evolves according to
\begin{align*}
R_3\quad:\quad
x_1\at{t}\approx x_2\at{t}\approx \ell\at{t}.
\end{align*} 
%
%
\subsubsection{Failure of the quasi-stationary approximation}\label{sect:312}
%
%
As illustrated in the right panel of Figure \ref{Fig:HyperbolicInstability}, the existence of the critical 
time $t_2$ can be understood  as follows. The quasi-stationary two-peaks approximation imposes the constraints
\begin{align*}
H^\prime\at{x_1}=H^\prime\at{x_2},\qquad x_1<x_2,\qquad x_*<x_2<x_{**},
\end{align*}
which define a smooth curve in the $\pair{x_1}{x_2}$-plane. This curve has the two branches
\begin{align*}
B_-:=\Big\{\bpair{X_-\at\si}{X_+\at\si}\;:\;\abs{\si}\leq\si_*\Big\}, \qquad
B_0:=\Big\{\bpair{X_0\at\si}{X_+\at\si}\;:\;\abs{\si}\leq\si_*\Big\}, \qquad
\end{align*}
which meet smoothly at the point $\pair{-x_*}{x_{**}}$. 
Due to the dynamical constraint the system is further confined to
\begin{align*}
G_\ell:=\Big\{\pair{x_1}{x_2}\;:\; m_1x_1+m_2x_2=\ell\Big\}\,,
\end{align*}
which is a straight line with slope $-m_1/m_2$. It can now easily be seen that in the quasi-stationary approximation for $t_0<t<t_1$ the state of the system corresponds to the unique intersection point of $B_-$ and $G_\ell$, which moves towards  $\pair{-x_*}{x_{**}}$ since $\ell$ is increasing in time. At time $t_1$, the system crosses this point  $\pair{-x_*}{x_{**}}$, and for sufficiently small times $t>t_1$ the state of the system is given by the unique intersection point of $B_0$ and $G_\ell$. At time $t_2$ however, the line $G_\ell$ becomes tangential to $B_0$
and the system can no longer follow the curve $B_0$ due to $\dot\ell>0$. Instead both peaks merge rapidly and the system jumps almost instantaneously to $\pair{\ell}{\ell}$, which is the only
intersection point of $G_\ell$ and the diagonal $x_1=x_2$. 
\par
The tangency condition for the merging reads
\begin{align*}
{m_1}{H^{\prime\prime}\at{x_2}}+
{m_2}{H^{\prime\prime}\at{x_1}}=0,
\end{align*}
where the left hand side is positive for all times $t_0<t<t_2$. Of course, if $m_2$ is very small, then the slope of $G_\ell$ is very negative and the tangency condition cannot be satisfied. In this case, there is a continuous
transition from the unstable-stable two-peaks configuration to the stable single-peak configuration since both peaks merge via $x_1=x_2=\ell$. It is then natural to define $t_2$ as the time at which this continuous merging takes place. 
\par
The above considerations also apply to the limiting case $m_1=1$ and $m_2=0$, 
provided that we accept that $x_2$ is undefined. In this case, the quasi-stationary approximation reads $x_1\at{t}=\ell\at{t}$ and the switching time $t_1$ is defined by $\ell\at{t_1}=-x_*$. At time $t_2$ with
$\ell\at{t_2}=+x_*$ the single peak leaves the unstable interval and in view of the above discussion  
it makes sense to interpret this event as (trivial) continuous merging.
\bigpar
We emphasize that our characterization of the unstable-stable evolution 
relies on condition (A3), that means on the concavity of the function $x_1\in\ccinterval{-x_*}{x_*}\mapsto x_2=X_+\at{H^\prime\at{x_1}}$. Otherwise it may happen at time $t_2$ that the system does not jump to the diagonal but instead to another unstable-stable configuration.
\par
We also mention that discontinuous merging of peaks implies that $\si$ jumps down to $H^\prime\at\ell$, and that the emerging single peak is stable due to $x_*<\ell<x_{**}$. Moreover, the energy of the system also jumps down. This is in line with the gradient flow structure of (\TPM), and holds because the energy decreases along the straight line that connects $\pair{x_1}{x_2}$ to $\pair{\ell}{\ell}$. In fact, with $E\at{\la}=m_1H\bat{\at{1-\la}x_1+\la\ell}+
m_2H\bat{\at{1-\la}{x}_2+\la\ell}$ we find $\dint{E}/\dint\la<0$ thanks to $\ell=m_1x_1+m_2x_2$, $x_1<x_2$
and $H^\prime\bat{\at{1-\la}x_1+\la\ell}<H^\prime\bat{\at{1-\la}x_2+\la\ell}$.
%
%
\subsubsection{Stability of unstable-stable two-peaks configurations}
%
%
We finally show that the quasi-stationary approximation is linearly stable until it ceases to exist. To this end, we consider a given solution to (\TPM) with $0<\tau\ll1$ and initial data as in \eqref{TPM.IV}, and denote by $\tilde{x}_1$ and $\tilde{x}_2$ the quasi-stationary approximation to (\TPM). This reads
\begin{align}
\label{TPM.QuasiStationaryApprox}
\tilde{x}_1\at{t}=X_{-/0}\bat{\tilde{\si}\at{t}},\qquad
\tilde{x}_2\at{t}=X_+\bat{\tilde{\si}\at{t}},\qquad
m_1\tilde{x}_1\at{t}+m_2\tilde{x}_2\at{t}=\ell\at{t},
\end{align}
where the definition of $\tilde{x}_1$ involves $X_-$  for $t_0\leq{t}\leq{\tilde{t}_1}$ but $X_0$ for $\tilde{t}_1\leq{t}\leq\tilde{t}_2$. Here, $\tilde{t}_2$ is the time at which the quasi-stationary approximation ceases to exist, and $\tilde{t}_1$ is defined by $\ell\at{\tilde{t}_1}=-m_1x_*+m_2x_{**}$ and denotes the switching time at which the quasi-stationary approximation enters the branch $B_0$.
\par%
Making the ansatz $x_i\at{t}=\tilde{x}_i\at{t}+\tau y_i\at{t}$, we easily derive the linearized model
\begin{align*}
\tau\dot{y}_1&=m_2\Bat{H^{\prime\prime}\at{\tilde{x}_2}y_2-H^{\prime\prime}\at{\tilde{x}_1}y_1}+g_1,\\
\tau\dot{y}_2&=m_1\Bat{H^{\prime\prime}\at{\tilde{x}_1}y_1-H^{\prime\prime}\at{\tilde{x}_2}y_2}+
g_2,\\
\si&=\tilde{\sigma}+\tau\at{m_1H^{\prime\prime}\at{\tilde{x}_1}y_1+
m_2H^{\prime\prime}\at{\tilde{x}_2}y_2+\dot\ell}
\end{align*}
with $g_i=\dot\ell-\dot{\tilde{x}}_i$.
By construction, we have $m_1y_1+m_2y_2=0$ and hence
\begin{align*}
\tau\dot{y}_i=-\zeta y_i+g_i,
\end{align*}
where  $\zeta=m_1H^{\prime\prime}\at{\tilde{x}_2}+m_2H^{\prime\prime}\at{\tilde{x}_1}$. Notice that $\zeta\at{t}$, $g_1\at{t}$ and $g_2\at{t}$ depend continuously on $t$ since the function $\ell$ is smooth. The Variations of Constants Formula now reveals that the quasi-stationary solution is dynamically stable as long as $\zeta>0$, that means as long as $t<\tilde{t}_2$. Moreover, a similar analysis reveals that the quasi-stationary single-peak solution $\tilde{x}_1=\tilde{x}_2=\ell$ is dynamically stable provided that $H^{\prime\prime}{\at\ell}>0$, which is satisfied for $t>\tilde{t}_2$ due to $\ell\at{\tilde{t}_2}\geq{x_{*}}$.
%
%
\subsection{Entropic effects: Widening and splitting of unstable peaks}\label{sect:32}
%
In the previous section we have seen that the two-peaks model (\TPM) allows a stable peak to enter the unstable region. Our goal in this section is to show that under some conditions the same is true for the original model (\FPM), but also that the entropic terms trigger new phenomena.
To point out the main idea we start with some heuristic arguments. Afterwards we 
derive and investigate a simplified model that allows to understand the key effects in a more formal way.
\bigpar%
Suppose that we are given a solution to (\FPM) that at time $t_0$ consists of two narrow stable peaks located at positions $-x_{**}<x_1\at{t_0}<-x_*$ and $x_*<x_2\at{t_0}<x_{**}$. Suppose also that the width of these peaks is sufficiently narrow such that  $\abs{\si\at{t_0}}<\si_*$. For small $\tau$ we then expect that the system
relaxes very fast to a meta-stable state of (\FPM) that meets the constraint $\ell\at{t_0}$. Without loss of generality we may hence assume that the initial peaks have width of order $\nu$ and that the initial data are
well prepared in the sense of \eqref{TPM.IV}
\par
For sufficiently small times $t>t_0$ we can expect that (\FPM) follows the two-peaks model (\TPM), that means
both peaks have width of order $\nu$ and are located at stable positions $x_1\in\oointerval{-x_{**}}{-x_*}$, $x_2\in\oointerval{x_*}{x_{**}}$ with $\dot{x}_1>0$ and $\dot{x}_2>0$  due to $\dot\ell>0$. At some time
$t_1$, however, the peak located at $x_1$ reaches the critical position $-x_*$, and this time can be estimated by
$\ell\at{t_1}\approx-m_1x_*+m_2x_{**}$. 
\par%
When the first peak has crossed $-x_*$, its width widens very quickly because the characteristics of the transport term now separate exponentially  with local rate $-H^{\prime\prime}{\at{x}}/\tau$. However, since the width of the first peak is initially exponentially small in $\tau$, it remains small for some times although it is surely much larger than $\nu$. Moreover, the second peak located at $x_2$ still has width of order $\nu$ as it remains confined to the stable region $x>x_*$.  Combining both arguments we conclude that the system can be approximated by (\TPM) even at times $t>t_1$ provided that $t-t_1$ is sufficiently small. The condition $\dot\ell>0$ then implies $\dot{x}_1>0$, $\dot{x}_2<0$ and $\dot\si<0$.
\par
\begin{figure}[ht!]%
\centering{%
\includegraphics[width=0.75\textwidth, draft=\figdraft]%
{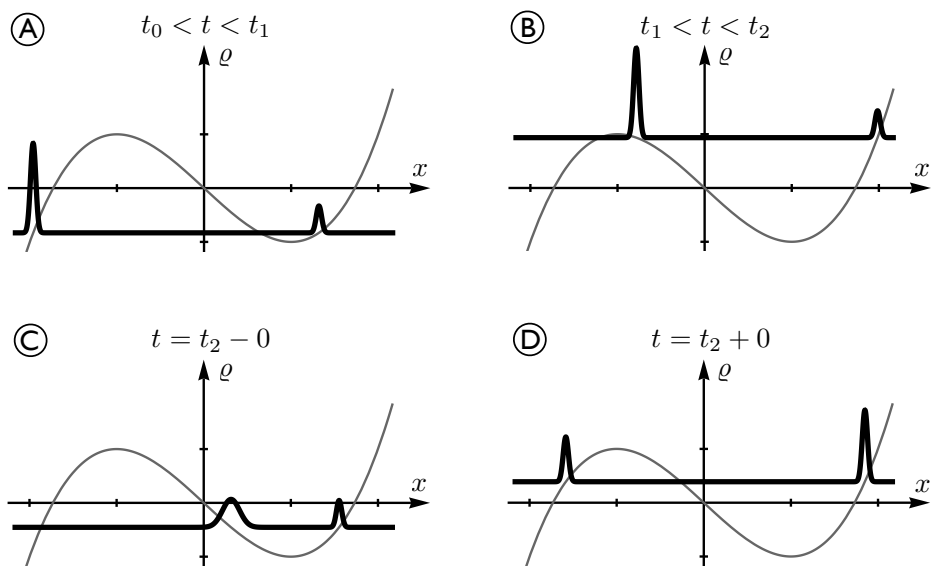}%
}%
\caption{%
Schematic representation of the shape of $\varrho$ at different times (for better illustration with respect to the graph of $H^\prime$, that means the plots show $\varrho\pair{x}{t}+\sigma\at{t}$ against $x$).
Due to the dynamical constraint, the peak located at $x_1$ switches from stable ($A$) to unstable ($B$) at
time $t_1$, and its width starts to grow exponentially. At time $t_2$, the width of the peak becomes suddenly of order $1$  ($C$) and the mass splitting problem describes that the system jumps almost instantaneously to a new stable-stable two-peaks configuration ($D$).
}%
\label{Fig:EntropicInstability}%
\end{figure}
The key question now is how long the first peak remains localized. Depending on the scaling parameter $a$ and the mass distribution between the peaks, it can happen that the first peak remains localized until both peaks merge continuously via $x_1\approx{}x_2\approx\ell\approx{}x_*$. In this case the first peak can in fact pass through the whole unstable region. It can also happen that the peak remains localized till the quasi-stationary two-peaks approximation ceases to exist. In this case both peaks merge instantaneously and discontinuously to form a single stable peak, but (\FPM) still behaves like (\TPM).
\par
There is, however, a third possible scenario, in which the width of the first peak becomes of order one before both
peaks can merge continuously or discontinuously. Below we will argue that if this happens at all, it happens
at a precise time $t_2$. More precisely, we show that there is a time $t_2$ at which the width of the first peak blows up instantaneously for small $\tau$. Some amount of the mass of this peak is then transported along characteristics to the left until it creates a new peak in the stable region $x<-x_*$. The remaining part, however, is transported towards
the other stable region $x>x_*$ to feed the second peak. Since the transport along characteristics if very fast, 
we expect that in the limit $\tau\to0$ the first peak splits and disappears instantaneously and that the system jumps
to another stable-stable configuration. After this jump the new peak in the stable region $x<-x_*$ has mass smaller than $m_1$, and this implies that the mass of the second peak is larger than $m_2$.
%
%
\subsubsection{A simplified model}
%
%
%
To analyze the relevant phenomena we study a simplified \emph{peak-widening model} (\PWM), which approximates the second peak at $x_2$ by a Dirac mass but keeps the probabilistic description for the first peak located at $x_1$. This gives rise to the equations
\begin{align}
\label{PWM-1}\tag{\PWM$_1$}
\tau\partial_t\varrho&=\partial_x\Bat{\nu^2\partial_x\varrho+\bat{H^\prime\at{x}-\si}\varrho}\,,\\
\label{PWM-2}\tag{\PWM$_2$}\tau\dot{x_2}&=\si-H^\prime\at{x_2}\,,\\
\label{PWM-3a}\tag{\PWM$_3$}
\ell&=m_1\int_\Rset{x}\varrho\dint{x}+m_2x_2\,,\\
\label{PWM-3b}\tag{\PWM$_3^{\,\prime}$}
\si&=m_1\int_\Rset{H^\prime\at{x}}\varrho\dint{x}+m_2x_2+\tau\dot\ell\,.
\end{align}
Here $m_1$ and $m_2=1-m_1$ are two constants that describe the mass distribution between the peaks, and as long as $x_2$ is confined to the stable region $x>x_*$ we can expect that each solution to (\PWM) defines via $m_1\varrho+m_2\delta_{x_2\at{t}}$ an (approximate) solution to the original model (\FPM).
\par%
In what follows we denote the width of the first peak by $w\at{t}$ and aim to derive
an asymptotic formula for the evolution of $w\at{t}$ that involves only the dynamical constraint $\ell\at{t}$. For simplicity we assume again
that the data at time $t_0$ are localized and well prepared in the sense of 
\begin{align*}
\varrho\pair{x}{t_0}\approx\delta_{x_1\at{t_0}}\at{x}, \qquad
w\at{t_0}\sim\nu,
\qquad -x_{**}<x_1\at{t_0}<-x_*,\quad
\qquad x_{*}<x_2\at{t_0}<x_{**}.
\end{align*}
A key ingredient to any asymptotic analysis of the widening phenomenon is to give an appropriate description of the position of the first peak. Our ansatz is to define $x_1\at{t}$ as solution of the characteristic ODE, that means we set
\begin{align}
\label{PW.Center1}
\tau\dot{x_1}&=\si-H^\prime\at{x_1}.
\end{align}
This ansatz has the following advantages. As long as the first peak is narrow, we have
\begin{align*}
\si\approx m_1H^\prime\at{x_1}+m_2H^\prime\at{x_2},\qquad
\ell\approx m_1x_1+m_2{x_2},\qquad
\end{align*}
and hence $x_1$ and $x_2$ evolve according to the two-peaks model (\TPM). Consequently, for $0<\tau\ll1$ we can describe $x_1$ and $x_2$ in terms of the quasi-stationary approximations of (TPM), whose dynamics is completely determined by the dynamical constraint $\ell\at{t}$ and the mass distribution between the peaks, see \eqref{TPM.QuasiStationaryApprox}. A further advantage of \eqref{PW.Center1} is that it gives rise to a quite simple evolution law for $w\at{t}$. 
%
%
\subsubsection{Formula for the width of the peak}
%
In order to analyze the growth of $w\at{t}$ we introduce the rescaling
\begin{align}
\label{PW.Scaling}
\varrho\pair{x}{t}=:\frac{1}{\nu\la\at{t}}R\Bpair{\frac{x-x_1\at{t}}{\nu\la\at{t}}}{\theta\at{t}}\,,
\end{align}
where both the spatial scaling factor $\la\at{t}$ and the rescaled time $\theta\at{t}$ will be identified below. We then have $w\at{t}=\nu\lambda\at{t}W\at{\theta\at{t}}$, where $W\at\theta$ is the width of $R$ at $\theta$.
\par
We readily verify that \eqref{PWM-1} and \eqref{PW.Center1} imply that the evolution of $R$ is governed by
\begin{align}
\label{PW.PDE1}
\tau\dot{\theta}\partial_{\theta}R=\partial_y\Bat{\frac{H^\prime\at{x_1-{\nu\la}y}-H^\prime\at{x_1}}{\nu\la}R}+\frac{1}{\la^2}\partial_y^2R+\frac{\tau\dot{\la}}{\la}\partial_y\at{yR}.
\end{align}
For the subsequent considerations we now assume that $x_1\at{t}$ is a given continuous function with $x_1\at{t_1}=-x_*$ and
\begin{align*}
-x_{**}<x_1\at{t}<-x_*
\quad\forall\quad t_0\leq{t}<{t_1}\,,
\qquad
{-x_*}<x_1\at{t}<x_*\quad\forall\quad t_1<t<t_3
\end{align*}
for some $t_3$. Heuristically, $t_3$ is the time at which both peaks merge (continuously or discontinuously) according to the two-peaks model (\TPM).
\par
As long as $w\at{t}$ is small -- but possibly much larger than $\nu$ --
we can expand the nonlinearity according to
\begin{align*}
H^\prime\at{x_1-{\nu\la}y}-H^\prime\at{x_1}=
H^{\prime\prime}\at{x_1}y+\DO{\nu\la{y}}=H^{\prime\prime}\at{x_1}y+\DO{w}.
\end{align*}
Neglecting the higher order terms and defining $\la$ and $\theta$ as solutions to
\begin{align*}
\tau\dot{\la}=-H^{\prime\prime}\at{x_1}
\la, \qquad \tau\dot{\theta}=\la^{-2},\qquad \la\at{t_0}=1,\qquad\theta\at{t_0}=0,
\end{align*}
the nonlinear PDE \eqref{PW.PDE1} transforms
into the heat equation
\begin{align}
\notag
\partial_\theta{R}=\partial_y^2R.
\end{align}
Consequently, for large $\theta$ the rescaled profile $R$ evolves in an almost self-similar manner, that means we can approximate
\begin{align}
R\pair{y}{\theta}\approx\frac{1}{\sqrt{4\pi\theta}}\exp\at{-\frac{y^2}{4\theta}}.
\label{PW.HeatEquation.Solution}
\end{align} 
Notice that this approximation implies $W\at\theta\sim\sqrt\theta$ for $\theta\gg1$ and holds as long as $W\at{t_0}$ is of order $1$. In order to characterize the width of the original peak it remains to understand how $\lambda$ and $\theta$ depend on $t$. A direct computation yields
\begin{align*}
\la\at{t}=\exp\at{-\frac{1}{\tau}\int_{t_0}^t f\at{s}\dint{s}},\qquad
\theta\at{t}=\frac{1}{\tau}\int_{t_0}^t
\exp\at{\frac{2}{\tau}\int_{t_0}^{\tilde{t}} f\at{s}\dint{s}}\dint{\tilde{t}},
\end{align*}
where $f$ abbreviates
\begin{align*}
f\at{t}:=H^{\prime\prime}\at{x_1\at{t}}.
\end{align*}
Due to $0<\tau\ll1$ and
\begin{align*}
f\at{t}>0\quad\forall\; 0<t<t_1,\qquad\qquad
\dot{\theta}\at{t}>0\quad\forall\;
0<t<t_3,
\end{align*}
we now infer that $\theta\at{t}\gg1$ for all 
$t_0<t<t_3$. This implies $w\at{t}=\nu\la\at{t}\sqrt{\theta\at{t}}$, and hence
\begin{align}
\label{PW.Width1}
w\at{t}^2&=\frac{1}{\tau}\int_{t_0}^t\exp\at{-\frac{2}{\tau}\int_{\tilde{t}}^tf\at{s}\dint{s}-\frac{2a}{\tau}}\dint{\tilde{t}},
\end{align}
where $a$ is the scaling parameter from \eqref{Slow:Scaling}.
\par
We have now identified an explicit formula for $w\at{t}$, which involves only the function $x_1\at{t}$. Notice that for small $\tau$ this function can be computed by the quasi-stationary  two-peaks approximation, whose evolution
is independent of $\tau$ and completely determined by $m_1$ and $\ell\at{t}$.
%
\subsubsection{Asymptotic description of the widening}
%
Formula \eqref{PW.Width1} can be restated as 
\begin{align}
\label{PW.Width2}
w\at{t}^2=\frac{1}{\tau}\exp\at{-\frac{2\phi\at{t}+2a}{\tau}}
\int_{t_0}^t\exp\at{\frac{2\phi\at{\tilde{t}}}{\tau}}\dint{\tilde{t}},
\end{align}
where the function 
\begin{align}
\label{PW.Width3}
\phi\at{t}:=\int_{t_1}^tf\at{s}\dint{s}
\end{align}
is non-positive on the interval $\oointerval{0}{t_3}$ since  $t_1$ is the only root of $f$ in this interval.
We now aim to show that this formula implies that the critical time $t_2$ is determined by the conditions
\begin{align}
\label{MSP.SplittingCondition}
\phi\at{t_2}+a=0,\qquad t_1<t_2 .
\end{align}
More precisely, if $\phi\at{t}$ attains the value $-a$ for some time $t_2$ with $t_1<t_2<t_3$, then 
in the limit $\tau\to0$ the width of the unstable peak explodes instantaneously at $t_2$.  If, however, $\phi\at{t}$ is larger than $-a$ for all times $t_1<t<t_3$, then the width of the peak is exponentially small in $\tau$ even at time $t_3$.
\par%
In order to prove these assertions we now derive rough estimates for $w\at{t}$. At first we consider $t_0<t\leq{t_1}$. In this case we have
$\phi\at{\tilde{t}}\leq\phi\at{t}$ for all $0\leq\tilde{t}\leq{t}$, so \eqref{PW.Width2} gives
\begin{align*}
w\at{t}^2\leq\frac{t}{\tau}\exp\at{-\frac{2a}{\tau}}=\frac{t\nu^2}{\tau}.
\end{align*}
Now let $t>t_1$ but suppose that $t$ is sufficiently small such that 
$\phi\at{t}+a=\delta>0$. Using $\phi\nat{\tilde{t}}\leq0$ we then we estimate
\begin{align*}
w\at{t}^2\leq\frac{t}{\tau}\exp\at{-\frac{2\delta}{\tau}},
\end{align*}
and find again that the width of the peak is exponentially small in $\tau$. 
Finally, we consider a time $t$ with $t_1<t<t_3$ and $\phi\at{t}+a=-\delta<0$.
To derive a lower bound for $w\at{t}$ we 
now employ the continuity of $f$ at $t_1$ as follows.  For $0<\eps<t-t_1$ and $t_1\leq\tilde{t}\leq{t_1+\eps}$ we estimate
\begin{align*}
\phi\at{\tilde{t}}\geq{}f\at{t_1}\at{\tilde{t}-t_1}-\Do{\eps}
\end{align*}
to find
\begin{align*}
\frac{1}{\tau}\int_{t_0}^{t}\exp\at{\frac{2\phi\at{\tilde{t}}}{\tau}}\dint{\tilde{t}}
&\geq
\frac{1}{\tau}\int_{t_1}^{t_1+\eps}\exp\at{\frac{2}{\tau}\Bat{f\at{t_1}\at{\tilde{t}-t_1}-\Do{\eps}}}\dint{\tilde{t}}
\\&\geq\frac{\exp\at{-\frac{2\Do{\eps}}{\tau}}}{2\abs{f\at{t_1}}}\at{1-\exp\at{\frac{2\eps{}f\at{t_1}}{\tau}}}
\xrightarrow{\eps\to0}\frac{1}{2\abs{f\at{t_1}}}>0.
\end{align*}
Combining this with \eqref{PW.Width2} gives
\begin{align*}
w\at{t}^2\geq\frac{1}{2\abs{f\at{t_1}}}\exp\at{\frac{2\delta}{\tau}},
\end{align*}
and we conclude that the width of the first peak is exponentially large in $1/\tau$.
\bigpar%
We finally emphasize that the equation for $w\at{t}$ can be simplified for $t\approx{t_2}$ as follows. Exploiting the continuity of $f$ at $t_2$ we find $\phi\at{t_2+\tau{s}}=-2a+f\at{t_2}s\tau+\Do{s\tau}$.
For each $s$ we therefore have
\begin{align*}
w^2\at{t_2+s\tau}=\bat{w^2\at{t_2}+\DO{\nu^2}}\exp\bat{-f\at{t_2}{s}+\Do{1}},
\end{align*}
where $\Do{1}$ means arbitrary small for small $\tau$. In particular, choosing $\tilde{t}_2$ with $\abs{\tilde{t}_2-t_2}=\Do{1}$ such that $w^2\at{\tilde{t}_2}^2=1$ we find
\begin{align}
\label{PW.AsymptoticWidth}
w^2\at{\tilde{t}_2+s\tau}=\bat{1+\Do{1}}\exp\bat{\beta{s}},\qquad
\beta=-f\at{t_2}>0.
\end{align}
%
%
\subsubsection{The mass splitting problem}
%
As explained above, at the critical time $t_2\approx\tilde{t}_2$ we expect that
the system undergoes a rapid transition from the unstable-stable configuration to a new stable-stable configuration.
In order to describe this transition, in particular, to predict the mass distribution between the emerging stable peaks, we propose to study a simplified \emph{mass-splitting model} (\MSM), which describes (\PWM) on the rescaled time scale $s=\at{t-\tilde{t}_2}/\tau$ in the limit $\nu\to0$. In other words, (\MSM) consists of the equations
\begin{align*}
\label{MSM-1}\tag{\MSM$_1$}
\partial_s\varrho&=\partial_x\bat{\bat{H^\prime\at{x}-\si}\varrho},\\
\label{MSM-2}\tag{\MSM$_2$}\frac{\dint}{\dint{s}}{x_2}&=\si-H^\prime\at{x_2},\\
\label{MSM-3}\tag{\MSM$_3$}\si&=m_1\int_\Rset{H^\prime\at{x}}\varrho\dint{x}+m_2x_2,
\end{align*}
which have no diffusion and satisfy the constraint via $\frac{\dint}{\dint{s}}\at{m_1x_1+m_2x_2}=0$.  Notice that this equation is the Wasserstein-gradient flow for the energy
\begin{align*}
\calE=\calH=m_1\int_\Rset H\at{x}\varrho\dint{x}+m_2H\at{x_2}.
\end{align*} 
In view of the above discussion we now impose asymptotic initial conditions at $s=-\infty$, which reflect that the mass splitting process starts in an unstable-stable configuration and that the unstable peak is a rescaled Gaussian due to the entropic randomness. To this end we identify $\tilde{t}_2=t_2$ and denote by $\tilde{x}_1$ and $\tilde{x}_2$ 
the quasi-stationary two-peaks approximation at  $t_2$, i.e. we have
\begin{align*}
\tilde{x}_1=X_0\at{\tilde\sigma},\qquad
\tilde{x}_2=X_+\at{\tilde\sigma},\qquad
m_1\tilde{x}_1+m_2\tilde{x}_2=\ell\at{t_2}.
\end{align*}
In accordance with \eqref{PW.Scaling}, \eqref{PW.HeatEquation.Solution} and \eqref{PW.AsymptoticWidth} we now require that
\begin{align}
\label{MSP.IV1}
\varrho\pair{x}{s}\quad\xrightarrow{s\to-\infty}\quad\frac{1}{2\sqrt{\pi}\exp\at{\beta{s}}}\exp\at{-{\frac{\at{x-\tilde{x}_1}^2}{4\exp\at{2\beta{s}}}}},\quad \beta=-H^{\prime\prime}\at{\tilde{x}_1}>0
\end{align}
weakly$*$ in the sense of probability measures, and that
\begin{align}
\label{MSP.IV2}
x_2\at{s}\quad\xrightarrow{s\to-\infty}\quad\tilde{x}_2.
\end{align}
The gradient flow structure of (\MSM) implies that each solution approaches a stable-stable configuration in the limit $s\to\infty$. This reads
\begin{align*}
x_2\at{s}\quad\xrightarrow{s\to+\infty}\quad\hat{x}_2,\qquad\qquad
\varrho\pair{x}{s}\quad\xrightarrow{s\to+\infty} \quad\at{m_1-m_{12}}\delta_{\hat{x}_1}+m_{12}\delta_{\hat{x}_2},
\end{align*}
where $-x_{**}<\hat{x}_1<-x_*$ and $x_*<\hat{x}_2$ denote the positions of the emerging stable peaks and 
$0\leq{}m_{12}\leq{}m_1$ is precisely the amount of mass that is transferred during the splitting
process from the unstable region into the stable region $x>x_*$. Of course, the asymptotic data at $s=+\infty$ must comply with
\begin{align*}
H^\prime\at{\hat{x}_1}=H^\prime\at{\hat{x}_2},\qquad
\at{m_1-m_{12}}\hat{x}_1+\at{m_2+m_{12}}\hat{x}_2=\ell=m_1\tilde{x}_1+
m_2\tilde{x}_2,
\end{align*}
but these conditions do not determine the three quantities $m_{12}$, $\hat{x}_1$, and $\hat{x}_2$ completely. This is not surprising and reflects that the amounts of mass that are transferred towards the stable regions depends crucially on the asymptotic shape of the unstable peak. In our case, however, this shape is a rescaled Gaussian
and therefore we expect that the data at $s=+\infty$ are uniquely determined by the data at $s=-\infty$. In particular, we conjecture that there is a \emph{unique} and \emph{continuous} mass transfer function $M$ such that
\begin{align}
\label{MSP.MDef}
m_{12}=M\pair{m_1}{\tilde{\si}}.
\end{align}
Both the existence and continuity of $M$ are not obvious because  the mass splitting problem involves two subtle limits. At first one has to show that the asymptotic condition \eqref{MSP.IV1} gives rise to a well-posed initial value problem at $s=-\infty$. Second, one has to guarantee that solutions do not drift as $s\to\infty$ along the 
connected one-parameter family of equilibrium solutions. A rigorous justification of the mass splitting function $M$ is beyond the scope of this paper and left for future research. 
\par
It is, however, possible 
to compute numerical approximations of $M$ using the method of characteristics. More precisely, restricting to $N$ characteristics $\xi_n\at{s}$  with $n=1\tdots{N}$, the mass 
splitting problem  (\MSM) can be approximated by
\begin{align}
\label{MSP.AppODE1}
\frac{\dint}{\dint{s}}{\xi}_n\at{s}=\si\at{s}-H^\prime\at{\xi_n\at{s}},\qquad
\frac{\dint}{\dint{s}}{x}_2\at{s}=\si\at{s}-H^\prime\at{x_2\at{s}},
\end{align}
with
\begin{align}
\label{MSP.AppODE2}
\si\at{s}=\frac{m_1}{N}\sum_{n=1}^NH^\prime\bat{\xi_n\at{s}}+m_2H^\prime\at{x_2\at{s}},\qquad
\ell=\frac{m_1}{N}\sum_{n=1}^N\xi_n\at{s}+m_2x_2\at{s}=\const.
\end{align}
Moreover, to mimic the asymptotic initial conditions \eqref{MSP.IV1} and \eqref{MSP.IV2}, we choose a small parameter $\eps$ and set
\begin{align*}
\xi_n\at{0}=\tilde{x}_1+\eps\mathrm{erf}^{-1}\at{n/N},\qquad x_2\at{0}=
\tilde{x}_2,
\end{align*}
with $\mathrm{erf}^{-1}$ being the inverse of $\mathrm{erf}\at{\xi}=\pi^{-1/2}\int_\infty^\xi\exp\at{-x^2}\dint{x}$. The resulting finite-dimensional initial value problem can be integrated numerically, for instance by means of the explicit Euler scheme, which satisfies the constraint $\ell=\const$ up to computational accuracy.  In the limit
\begin{align*}
x_2\at{s}\quad\xrightarrow{s\to\infty}\quad\hat{x}_2,\qquad\xi_n\at{s}\quad\xrightarrow{s\to\infty}\quad\left\{
\begin{array}{lcl}%
\hat{x}_1&\text{for}& 1\leq{n}<N_{12},\\
\hat{x}_2&\text{for}&N_{12}<n\leq{N}.
\end{array}%
\right.
\end{align*}
The critical index $N_{12}$ finally determines $m_{12}$ via $m_{12}=m_1N_{12}/N$.
\par
\begin{figure}[t!]
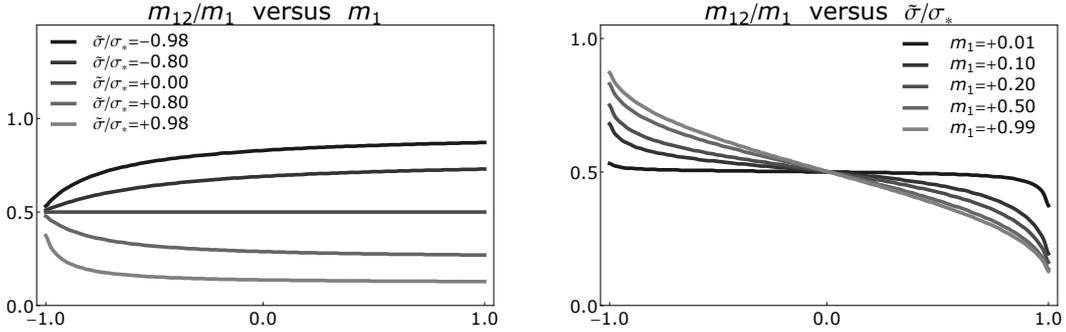
%
\centering{%
\includegraphics[width=0.4\textwidth, draft=\figdraft]%
{\figfile{msp_sim_m1}}%
\hspace{0.05\textwidth}
\includegraphics[width=0.4\textwidth, draft=\figdraft]%
{\figfile{msp_sim_si}}%
}%
\caption{%
Numerically computed values of $m_{12}$ 
for potential \eqref{ExamplePot} and several values of $m_1$ and $\tilde\si$.
}%
\label{Fig:msp_sim}%
\end{figure}

The strategy for the numerical computation of $M$ is therefore as follows. We sample the two-dimensional parameter
space of (\MSM) and solve for each choice of the parameters the ODE \eqref{MSP.AppODE1} with $\si$ as in \eqref{MSP.AppODE2} on a sufficiently large time interval. The results are illustrated in Figure \ref{Fig:msp_sim} for the potential \eqref{ExamplePot}, where $0\leq{m_1}\leq1$ and $\tilde\si=H^\prime\at{\tilde{x}_i}$ are, as above, regarded as the two independent parameters.
%
%
\subsection{Limit dynamics: Switching, splitting, and merging events}\label{sect:33}
%
%
Combining the arguments from the previous sections we are now able to characterize the
slow reaction limit of (\FPM) in terms of single-peak and two-peaks configurations. To describe their evolution 
we consider piecewise constant mass functions $m_1\at{t}$ and $m_2\at{t}$ along with
piecewise continuous position functions $x_1(t)$ and $x_2\at{t}$, where we
allow $x_1$ and $x_2$ to be undefined on intervals with $m_1=0$ and $m_2 =0$, respectively.
These functions are coupled by the constraints 
\begin{align}
\label{LM.Constraints}
H^\prime\at{x_1}=H^\prime\at{x_2}=\si\,,\qquad 
\ell=m_1x_1+m_2x_2\,,\qquad m_1+m_2=1
\end{align}
with
\begin{align*}
x_1\leq{}x_2\,,\qquad 0\leq
m_1,\,m_2\leq1\,,
\end{align*}
and define the limit energy via
\begin{align*}
E:=m_1H\at{x_1}+m_2H\at{x_2}.
\end{align*}
To formulate the limit model we first specialize to strictly increasing dynamical constraints and discuss possible generalizations afterwards.
%
%
\subsubsection{Limit model for strictly increasing constraints}
%
%
Sections~\ref{sect:31} and \ref{sect:32} reveal that the slow reaction limit for constraints with \eqref{Intro.Constraint} can be described as follows. The initial dynamics is governed by the quasi-stationary evolution of a single stable peak within the stable region $x<-x_*$ and hence it seems natural to define $m_1\at{t}=1$, which implies $x_1\at{t}=\ell\at{t}$ and renders $x_2$ to be undefined.  At the first switching time, which is defined by $x_1=\ell=-x_*$, the single peaks becomes unstable and its width starts to widen. It may happen that the motion of the unstable peak is faster than the widening, and then the peak remains localized until it leaves the unstable region via $\ell=x_*$. In this case we find a Type-I phase transition as there is no splitting of unstable peaks; notice that this happens if $\dot{\ell}$ is sufficiently large or if the scaling parameter $a$ is sufficiently large.
\par
For Type-II transitions, however, the width of the peak becomes eventually large within the unstable interval. This means there exists a first splitting time with $-x_*<\ell<x_*$ at which the system jumps from the unstable single-peak configuration to a stable-stable two-peaks configuration, where the masses and the positions of the emerging stable peaks is determined by the mass splitting problem. In particular, after the first splitting event the positions $x_1$, $x_2$ and the masses $m_1$, $m_2$ are well defined, and both stable peaks move according to the dynamical constraint until the first peak located at $x_1$ becomes unstable at the next switching time defined by $x_1=x_*$. Afterwards we must carefully investigate the unstable-stable evolution in order to decide whether there is a further splitting of the unstable peak or whether both peaks finally merge continuously or discontinuously.
\par
In summary, the limit dynamics can be described as illustrated
in Figure \ref{Fig:flow_chart}, that means intervals of quasi-stationary transport are interrupted 
by singular times corresponding to the following types of events:
\begin{list}{$\bullet$}
{ %
\setlength{\partopsep}{0pt}     
\setlength{\leftmargin}{6\parindent}
\setlength{\rightmargin}{\parindent}
\setlength{\labelwidth}{3\parindent}
\setlength{\itemindent}{0\parindent}
\setlength{\labelsep}{1\parindent} } 
\item[\emph{Switching}:] 
The peak at $x_1$ enters the unstable region.     
\item[\emph{Splitting}:] 
The unstable peak at $x_1$ splits and the system jumps to a new stable-stable configuration with decreased mass $m_1$ and increased mass $m_2$.
\item[\emph{Merging}:] 
The peaks in an unstable-stable configuration merge either continuously (with $x_1=x_*=x_2$) or discontinuously (with $x_1<x_*<x_2$), or there is only a single peak with $m_1=1$ that leaves  the unstable region (with $x_1=x_*$).
\end{list}%
More precisely, each limit trajectory comprises $K+1$ switching times,
$K$ splitting times, and a final merging time
\begin{align*}
0<t_{\evsw,0}<t_{\evsp,1}<t_{\evsw,2}<\tdots<t_{\evsp,K}<t_{\evsw,K}<t_{\evme}<t_*\,.
\end{align*}
where we have $K=0$ and $K\geq1$ for Type-I and Type-II transitions, respectively. Here $t_*<\infty$ is defined by $\ell\at{t_*}=x_{**}$, so proper two-peaks configurations can only exist for $t<t_*$. 
\par
Notice that an infinite number of switching and splitting events is not possible because the splitting condition \eqref{MSP.SplittingCondition} implies a lower bound for $t_{\evsp,k}-t_{\evsw,k-1}$ via
\begin{align}
\notag
a=\int_{t_{\evsw,k-1}}^{t_{\evsp,{k}}}\abs{H^{\prime\prime}\at{x_1\at{t}}}\dint{t}\leq C\at{t_{\evsp,k}-t_{\evsw,k}}
\end{align} 
with $C=\sup_{\abs{x}\leq{x_*}} \abs{H^ {\prime\prime}\at{x}}$. Notice also that we must truly book keep the switching events since  \eqref{MSP.SplittingCondition} doe not determine $t_{\evsp,k}$ but only
$t_{\evsp,k}-t_{\evsw,k-1}$. 
\bigpar
We now describe the flowchart from Figure \ref{Fig:flow_chart} in greater detail. 
\mhparagraph{Intervals of transport} Between consecutive singular times, and likewise initially for $0<t<t_{\evsw,0}$, the dynamics of $x_1$ and $x_2$ is governed by a rate-independent system of non-autonomous ODEs. More precisely, differentiating \eqref{LM.Constraints} with respect to time yields
\begin{align}
\label{LM.ODEs1}
\dot{m}_1=\dot{m}_2=0
,\qquad
\dot{x}_1=\frac{H^{\prime\prime}\at{x_2}}{Z\quadruple{m_1}{m_2}{x_1}{x_2}}\dot\ell
,\qquad
\dot{x}_2=\frac{H^{\prime\prime}\at{x_1}}{Z\quadruple{m_1}{m_2}{x_1}{x_2}}\dot\ell
,%
\end{align}
where we recall that $x_1$ and $x_2$ remain undefined for $m_1=0$ and $m_2=0$, respectively. These ODEs imply
\begin{align}
\label{LM.ODEs2}
\dot{E}=\si\dot\ell, \qquad \dot\si=
\frac{H^{\prime\prime}\at{x_1}H^{\prime\prime}\at{x_2}}{Z\quadruple{m_1}{m_2}{x_1}{x_2}}\dot\ell,
\end{align}
with
\begin{align*}
Z\quadruple{m_1}{m_2}{x_1}{x_2}=m_1H^{\prime\prime}\at{x_2}+
m_2H^{\prime\prime}\at{x_1}.
\end{align*}
Since the dynamical constraint $\ell\at{t}$ is given, the initial value problem 
to \eqref{LM.ODEs1} can -- at least in principle -- be integrated.
\par
\mhparagraph{Jump conditions at singular times}
At switching times all functions are continuous, that means we have 
\begin{align*}
\jump{m_i}\at{t_{\evsw,k}}=\jump{x_i}\at{t_{\evsw,k}}=
\jump{\si}\at{t_{\evsw,k}}=\jump{E}\at{t_{\evsw,k}}=0\,,\qquad i=1,2\,,
\end{align*}
where $\jump{f}\at{t}:= f\at{t+0}-f\at{t-0}$ with $f\at{t\pm0}=\lim_{\delta\searrow0}f\at{t\pm\delta}$ denotes the jump of the function $f$ at time $t$.
\par%
At a splitting time all quantities do jump, where the jump heights are determined by the mass splitting problem (\MSM). In particular, 
the gradient flow structure of (\MSM) ensures that 
\begin{align*}
\jump{E}\at{t_{\evsp,k}}<0,
\end{align*}
and using the mass transfer function $M$ from \eqref{MSP.MDef} we find
\begin{align*}
\jump{m_2}\at{t_{\evsp,k}}=-\jump{m_1}\at{t_{\evsp,k}}=
M\bpair{m_1\at{t_{\evsp,k}-0}}{\si\at{t_{\evsp,k}-0}}.
\end{align*}
The jumps of $x_1$, $x_2$ and $\si$ are then determined by the algebraic constraints \eqref{LM.Constraints}, and can be used to reinitialize the ODEs \eqref{LM.ODEs1} and \eqref{LM.ODEs2}. 
\par
Finally, at the merging time we have
\begin{align*}
x_1\at{t_{\evme}+0}=x_2\at{t_{\evme}+0}=\ell\at{t_{\evme}}\,,\qquad
\si\at{t_{\evme}+0}=H^\prime\bat{\ell\at{t_{\evme}}}\,,
\end{align*}
and the concavity condition (A3), see the discussion at the end of Section \ref{sect:312}, provides
\begin{align*}
\jump{E}\at{t_{\evme}}\leq0
\end{align*}
with strict inequality for discontinuous merging. After the merging 
the precise values masses $m_1$ and $m_2$ are actually undefined, but it seems natural to set
\begin{align*}
m_1\at{t_{\evme}+0}=0,\qquad m_2\at{t_{\evme}+0}=1\,.
\end{align*}

\mhparagraph{Determining the next singular time}
After the $k^\text{th}$ splitting time (we set $t_{\evsp,0}:=0$ to describe the initial evolution),  the system is in a stable-stable configuration and the subsequent
switching time $t_{\evsw,k}$ is the smallest time larger than $t_{\evsp,k}$ such that $\si\at{t_{\evsw,k}}=\si_*$. This implies $x_1\at{t_{\evsw,k}}=-x_{*}$, $x_2\at{t_{\evsw,k}}=x_{**}$, and
\begin{align*}
\ell\at{t_{\evsw,k}}=-m_1x_{*}+m_2x_{**},
\end{align*}
which determines $t_{\evsw,k}$ uniquely since $m_1$ and $m_2$ are known. 
\par
After the switching event at $t_{\evsw,k-1}$ we have to decide whether the next singular time corresponds to splitting or merging according to the conditions
\begin{align*}
\text{\emph{discontinuous merging}:}&\qquad
Z\quadruple{m_1}{m_2}{x_1\at{t_{\evme}}}{x_2\at{t_{\evme}}}=0,\\
\text{\emph{continuous merging}:}&\qquad
x_1\at{t_\evme}=x_2\at{t_\evme}=\ell\at{t_\evme},\\
\text{\emph{splitting:}}&\qquad\int_{t_{\evsw,k-1}}^{t_{\evsp,k}}H^{\prime\prime}\bat{x_1\at{t}}\dint{t}+a=0.
\end{align*}
To discuss this decision in a simple case, let us suppose that $\dot\ell$ is constant in time. Due to $\dot{\si}<0$ we can replace $t$ by $\sigma$, and direct computations  provide the following nonlinear equations
\begin{align*}
\text{\emph{discontinuous merging}:}&\qquad
Z\quadruple{m_1}{m_2}{X_0\at{\si_{\evme}}}{X_+\at{\si_{\evme}}}=0,\\
\text{\emph{continuous merging}:}&\qquad
\si_{\evme}=-\si_*,\\
\text{\emph{splitting:}}&\qquad\int_{\si_{\evsp, k}}^{\si_*}
\frac{Z\quadruple{m_1}{m_2}{X_0\at{\tilde\si}}{X_+\at{\tilde\si}}}{H^{\prime\prime}\at{X_+\at{\tilde\si}}}\dint{\tilde\si}=a\dot\ell.
\end{align*}
In particular, both merging conditions are rate independent, whereas the switching condition depends on $\dot\ell$ and cannot be satisfied for large $\dot\ell$. This is not surprising because the merging conditions are completely
determined by the quasi-stationary two-peak approximation, whereas splitting only happens if the peak widening due to the separation of characteristics is faster than the transport due to the dynamical constraint. 
\par
Instead of  solving \eqref{LM.ODEs1}, we can now vary $\si$, starting from $\si_*$ and moving towards $-\si_*$, and check which conditions is satisfied at first. In any case, the time $t_\ev$ of the next singular event can be computed by
\begin{align*}
t_{\ev}-t_{\evsw,k}=\int_{\evsw,k-1}^{t_{\ev}}\dint{t}=\tfrac{1}{\dot\ell}%
\int_{\si_{\ev}}^{\si_*}
\frac{Z\quadruple{m_1}{m_2}{X_0\at{\si}}{X_+\at{\si}}}{H^{\prime\prime}\at{X_0\at{\si}}H^{\prime\prime}\at{X_+\at{\si}}}\dint{\si}\,.
\end{align*}
%
%
%
\subsubsection{Limit model for non-monotone constraints}\label{sect:332}
%
We finally derive an alternative description of the slow reaction limit that covers arbitrary
dynamical constraints $\ell$. Decreasing constraints can produce
stable-unstable configurations with $x_1<-x_*$ and $-x_*<x_2<x_*$, which in turn can split or
merge. However, these effects can easily be described by adopting the formulas from sections \ref{sect:31} and \ref{sect:32}. The truly new effect is that non-monotone constraints
can trigger \emph{inverse switching events}, that means, for instance, an
unstable-stable configuration can become stable-stable via $x_1=-x_*$ 
with $\dot{x}_1<0$.
\par
For general constraints, it is convenient to describe the slow-reaction limit in terms of the following variables:
$\at{i}$ the dynamical multiplier $\si$, $\at{ii}$ an internal variable $\phi\in\ccinterval{-a}{0}$ as introduced in Section \ref{sect:32} to control the width of an unstable peak, and $\at{iii}$ 
three nonnegative masses $m_-$, $m_0$, and $m_+$ with
\begin{align*}
\sum_{i\in\{-,0,+\}}m_i=1\,,\qquad \prod_{i\in\{-,0,+\}}m_i=0
\end{align*}
to describe the several types of single-peak and two-peaks configurations. Using these parameters, 
the peak positions are dependent variables which satisfy $x_i=X_i\at{x}$ for $i\in\{-,0,+\}$ as long as $X_i$ is well-defined.
\begin{table}[ht!]%
\centering{%
\begin{tabular}{ccccc}%
configuration&masses $m_i$&range for $\si$&
range for $\phi$&further constraint
\\\hline\\%
$\calS_{-}$&$m_-=1$&$-\infty<\si< + \si_*$&$\phi=0$&
\\%
$\calS_{+}$&$m_+=1$&$-\si_*<\si< +\infty$&$\phi=0$&
\\%
$\calS_{0}$&$m_0=1$&$-\si_*<\si< +\si_*$&$-a<\phi<0$&
\medskip\\%
$\calT_{-+}$&$m_0=0$&$-\si_*<\si< + \si_*$&$\phi=0$&
\\%
$\calT_{-0}$&$m_+=0$&$-\si_*<\si< + \si_*$&$-a<\phi<0$&$Z_{-0}>0$
\\%
$\calT_{0+}$&$m_-=0$&$-\si_*<\si< + \si_*$&$-a<\phi<0$&$Z_{0+}>0$
\end{tabular}%
}%
\caption {List of all possible single-peak configurations $\calS_i$ and two-peaks configurations $\calT_{ij}$.}%
\label{Tbl:Sets}%
\end{table}%
\par%
At each non-singular time, the system is confined to one of the sets defined in Table \ref{Tbl:Sets}, where
$Z_{-0}$ and $Z_{0+}$ are abbreviations for
\begin{align*}
Z_{-0}:=m_-A_0\at{\si}+m_0A_-\at{\si}\,,\qquad
Z_{0+}:=m_0A_+\at{\si}+m_+A_0\at{\si}\,,
\end{align*}
and the functions $A_i$ are defined by
\begin{align*}
A_i\at\si:=H^{\prime\prime}\bat{X_i\at{\si}},\qquad i\in\{-,0,+\}\,.
\end{align*}%
Within each of these sets, the peaks are transported according to the quasi-stationary
two-peaks approximation and the widening of unstable peaks is governed by $\phi$ as described in \eqref{PW.Width2} and \eqref{PW.Width3}. This reads
\begin{align}
\label{LM:RegularDynamics}
\dot{m_i}=0,\,\qquad\dot{\si}=\dot\ell\!\!\!\!\!\sum_{i\in\{-,0,+\}}\!\!\!\!\!m_iA_i\at{\si},\,\qquad
\dot\phi=\chi_{\{m_0\neq0\}}A_0\at\si\,.
\end{align}
\begin{table}[ht!]%
\centering{%
\begin{tabular}{llll}%
singular event&condition&$\quad$&possible jumps
\\\hline\\%
switching&$\si=+\si_*$&&
$\calS_-\mapsto\calS_0,\quad \calT_{-+}\mapsto\calT_{0+}$
\medskip\\%
&$\si=-\si_*$&&
$\calS_+\mapsto\calS_0,\quad\calT_{-+}\mapsto\calT_{-0}$
\medskip\\%
inverse switching&$\si=+\si_*$&&
$\calS_0\mapsto\calS_-,\quad\calT_{0+}\mapsto\calT_{-+}$
\medskip\\%
&$\si=-\si_*$&&
$\calS_0\mapsto\calS_+,\quad\calT_{-0}\mapsto\calT_{-+}$
\medskip\\%
splitting&$\phi=-a$&&$\calS_{0}\mapsto\calT_{-+},\quad\calT_{-0}\mapsto\calT_{-+},\quad\calT_{0+}\mapsto\calT_{-+}$
\medskip\\%
discontinuous merging&$Z_{-0}=0$&&$\calT_{-0}\mapsto\calS_{-}$
\medskip\\%
&$Z_{0+}=0$&&$\calT_{0+}\mapsto\calS_{+}$
\medskip\\%
continuous merging&$\si=+\si_*$&&$\calT_{-0}\mapsto\calS_{-}$
\medskip\\%
&$\si=-\si_*$&&$\calT_{0+}\mapsto\calS_{+}$%
\end{tabular}%
}%
\caption {List of all possible singular events.}%
\label{Tbl:Events}%
\end{table}%
\begin{table}[ht!]%
\centering{%
\begin{tabular}{lll}%
singular event&subcase&behavior of variables%
\\\hline\\%
switching&$\si=+\si_*$&$m_- \leftrightarrow m_0$%
\medskip\\%
&$\si=-\si_*$&$m_0 \leftrightarrow m_+$%
\medskip\\%
inverse switching&$\si=+\si_*$&$m_- \leftrightarrow m_0,\quad \phi\to0$%
\medskip\\%
& $\si=-\si_*$& $m_0 \leftrightarrow m_+,\quad \phi\to0$%
\medskip\\%
splitting&&
$m_0\to0,\quad \jump{m_-}\geq0, \quad \jump{m_+}\ge0, \quad \phi\to0,\quad\jump{\si}\gtreqqless0$%
\medskip\\%
merging&into $\calS_+$&$m_+\to m_0+m_+,\quad m_0\to0,\quad \phi\to0,\quad
\jump{\si}\le0$%
\medskip\\%
&into $\calS_-$& $m_-\to m_-+m_0,\quad m_0\to0,\quad \phi\to0,\quad
\jump{\si}\ge0$%
\end{tabular}%
}%
\caption {Update rules for variables $m_i$, $\si$, and $\phi$ at singular events.}%
\label{Tbl:Jumps}%
\end{table}%
On the other hand, singular events happen when the system reaches the boundary of either $\calS_i$ or $\calT_{ij}$. More precisely, examing the conditions for switching, inverse switching, splitting, and merging we arrive at
list from Table \ref{Tbl:Events}, where we now interpret trivial continuous merging as inverse switching of single-peak configurations. The corresponding jump and update rules for the variables $m_i$, $\si$, and $\phi$
are summarized in Table \ref{Tbl:Jumps}. 
\par
The slow reaction dynamics of two-peaks initial data can now be integrated iteratively by $\at{i}$ 
moving via \eqref{LM:RegularDynamics} along either one of the sets $\calS_i$ and $\calT_{ij}$, and $\at{ii}$ 
jumping to another set when reaching the boundary. Notice that at both splitting and discontinuous merging events the system jumps to inner points and that splitting events require to solve the mass splitting problem.
%
%
%
%

\section*{Acknowledgements}%
The authors are grateful to Wolfgang Dreyer, Clemens Guhlke, and Alexander Mielke for
stimulating discussions. This work was supported by the EPSRC Science and Innovation award to the Oxford Centre for
Nonlinear PDE (EP/E035027/1).
%
%
%
%
\providecommand{\bysame}{\leavevmode\hbox to3em{\hrulefill}\thinspace}
\providecommand{\MR}{\relax\ifhmode\unskip\space\fi MR }
\providecommand{\MRhref}[2]{%
  \href{http://www.ams.org/mathscinet-getitem?mr=#1}{#2}
}
\providecommand{\href}[2]{#2}

\end{document}